\documentclass[12pt]{article}

\usepackage[intlimits]{amsmath} 
\usepackage{amsthm,mathrsfs}    
\usepackage[comma,authoryear]{natbib}
\usepackage{url}
\numberwithin{equation}{section}
\usepackage{amssymb,amsmath}
\usepackage{subfigure}
\usepackage{bm}
\usepackage{booktabs} 
\usepackage{graphicx,graphics,epsfig}
\usepackage{float}
\usepackage{wrapfig}

\theoremstyle{plain}
\newtheorem{thm}{Theorem}

\newtheorem{lemma}{Lemma}[thm]
\newtheorem{coro}{Corollary}[thm]

\newcommand{\bthm}{\begin{thm}}
\newcommand{\ethm}{\end{thm}}

\newcommand{\bpf}{\begin{proof}}
\newcommand{\epf}{\end{proof}}

\theoremstyle{definition}

\newtheorem{rem}{Remark}

\newtheorem{exmp}{Example}

\usepackage[a4paper,text={16.5cm,25cm},centering]{geometry}
\usepackage[compact]{titlesec}
\newcommand{\bib}{\bibliography{ref-bib}\bibliographystyle{Chicago}}
\usepackage{deep-macro}

\begin{document}
\pagestyle{empty}
\begin{center}
{\Large {\bf LP Approach to Statistical Modeling}}
\\[.2in]
Subhadeep Mukhopadhyay$^1$,  Emanuel Parzen$^2$, \\
$^1$Temple University, Philadelphia, PA, USA\\
$^2$Texas A\&M University, College Station, TX, USA\\[.35in]
{\bf ABSTRACT}\\
\end{center}
We present an approach to statistical data modeling and exploratory data analysis called `LP Statistical Data Science.'  It aims to generalize and unify traditional and novel statistical measures, methods, and exploratory tools. This article outlines fundamental concepts along with real-data examples to illustrate how the `LP Statistical Algorithm' can systematically tackle different varieties of data \emph{types}, data \emph{patterns}, and data \emph{structures} under a coherent theoretical framework. A fundamental role is played by specially designed orthonormal basis of a random variable $X$ for linear (Hilbert space theory) representation of a general function of $X$, such as $\Ex[Y \mid X]$.

\vspace*{.45in}


\vspace*{.5in}
\thispagestyle{empty}

\thispagestyle{empty}
\medskip\hrule height 1.4pt
{\small
\tableofcontents
}

\newpage
\pagestyle{plain}
\setcounter{page}{1}
\setcounter{section}{-1}
\section{Goals}
There has been explosive growth of statistical algorithms in recent times, due to greater availability of data with mixed types (discrete or continuous), patterns (simple linear to complex nonlinear) and structures (univariate to multivariate). Growth in statistical software technology contributed substantially to simplify the implementation workloads of these various methods. While execution of an algorithm getting easier by pushing the button, the theory and practice of statistical learning methods are becoming more complicated day by day with the availability of thousands of `isolated' statistical ideas and softwares. We seek a more systematic and automatic approach to develop interpretable algorithms that permits easy way to establish relationship among various statistical methods (thus simplifies teaching and facilitate practice). An attempt is made in this paper to describe how this can be done via a new theory `LP Statistical Data Science'.
\vskip.35em
Our theory provides unified representation of statistical methods based on few fundamental functions (defined for both discrete and continuous random variables): Quantile function, custom-constructed LP orthonormal score functions, LP Skew density function, LP copula density function and LP comparison density function. We apply these tools to solve broad range of statistical problems like goodness-of-fit, density estimation, correlation learning, conditional mean and quantile modeling, copula estimation, categorical data modeling, correspondence analysis etc., in a coherent manner.
\vskip.35em
We hope LP statistical learning theory might give us the clue for designing `United Statistical Algorithm'. We outline core mathematical concepts of the LP approach to data analysis, design applicable algorithms and illustrate using concrete real data examples.
\section{LP Score Function}
Our approach to nonparametric modeling is based on the representation of a general function $h(X)$ in an orthonormal basis of functions of $X$ (more precisely, functions of $\Fm(x;X)$, the mid-distribution function of $X$), denoted by $T_j(x;X)$, called LP Score functions. LP unit score functions are orthonormal on unit interval defined as $S_j(u;X)=T_j(Q(u;X);X)$, $0<u<1$, where $Q(u;X)$ is quantile function of $X$.
\subsection{Definition and Construction}
Define Mid-distribution function of a random variable $X$ as
\beq
\Fm(x;X)\,=\,F(x;X)-.5p(x;X),\,\, p(x;X)\,=\,\Pr[X=x], \,\,F(x;X)\,=\,\Pr[X \leq x].
\eeq
For a mid-distribution $\Fm(x;X)$ with quantile function $Q(u;X), 0<u<1$, we construct $T_j(x;X)$, $j=0,1,\ldots$ orthonormal LP score functions by Gram Schmidt orthonormalization of powers of $T_1(x;X)$. Define $T_0(x;X)=1$ and
\beq T_1(x;X) \,=\,\dfrac{\Fm(x;X) - .5}{\si(\Fm(X;X))}. \eeq
\cite{parzen04b} showed that $\Ex[\Fm(X;X)]=.5$ and $\var[\Fm(X;X)]=(1-\sum_x p^3(x;X))/12$. Note that for $X$ continuous we have $T_1(x;X)=\sqrt{12}(F(x;X)-.5)$.
LP unit score functions $S_j(u;X)$ are orthonormal in  $L^2(F)$, the Hilbert space of square integrable functions on $0<u<1$, with respect to the measure $F$ (either discrete or continuous). When X is continuous, we recover the Legendre polynomials
\beq
S_j(u;X)\,=\,\Leg_j(u),\,\,\,\, T_j(X;X)\,=\,\Leg_j\big[F(X;X)\big].
\eeq
Formulas of orthonormal shifted Legendre polynomials on $[0,1]$ are given in Appendix A1.

\subsection{Sufficient Statistics Interpretation}

Models that are nonlinear in $X$ can be constructed in many applications as linear models of the score vector $TX$ with components $T_j(X;X)$, which we interpret as ``sufficient statistic'' for $X$. This approach is implemented for time series analysis in our theory of LPTime \citep*{D12e}.

\begin{figure}[ht]
 \centering
 \vspace{-.15em}
 \subfigure[]{
  \includegraphics[height=.44\textheight,width=\textwidth,keepaspectratio,trim=1cm 1cm 1cm 1cm]{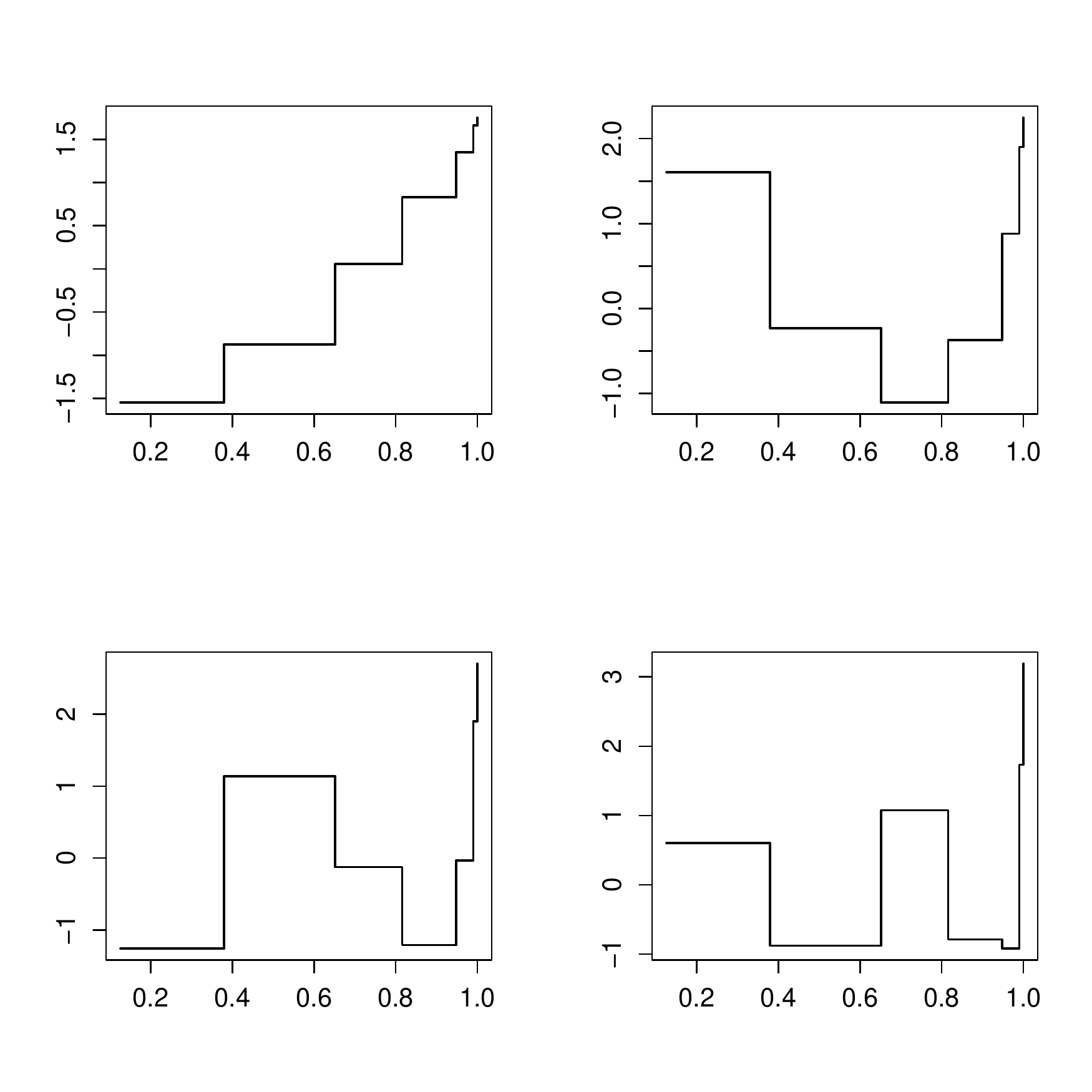}
   }
 \subfigure[]{
  \includegraphics[height=.44\textheight,width=\textwidth,keepaspectratio,trim=1cm 1cm 1cm 1cm]{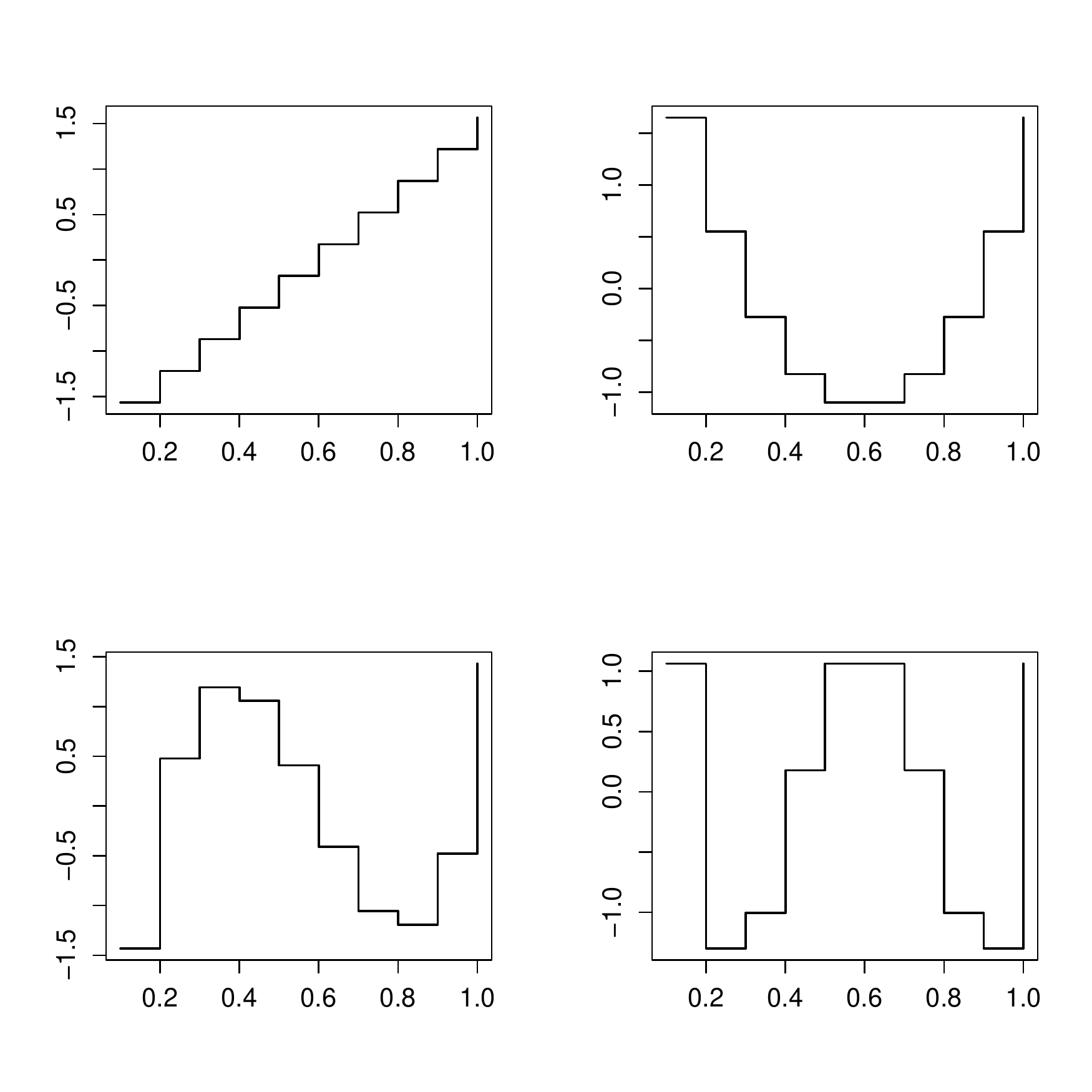}
   }
 \caption{Shapes of LP-unit score functions for (a) poisson distribution (with $\la=2$), and (b) discrete uniform distribution (which generates discrete Legendre polynomials) based on $n=250$ simulated samples.  }
  \label{fig:score}
\end{figure}

\begin{figure}[tbh!]
 \centering
 \includegraphics[height=.42\textheight,width=\textwidth,keepaspectratio,trim=1.5cm 1cm 1.5cm 1.5cm]{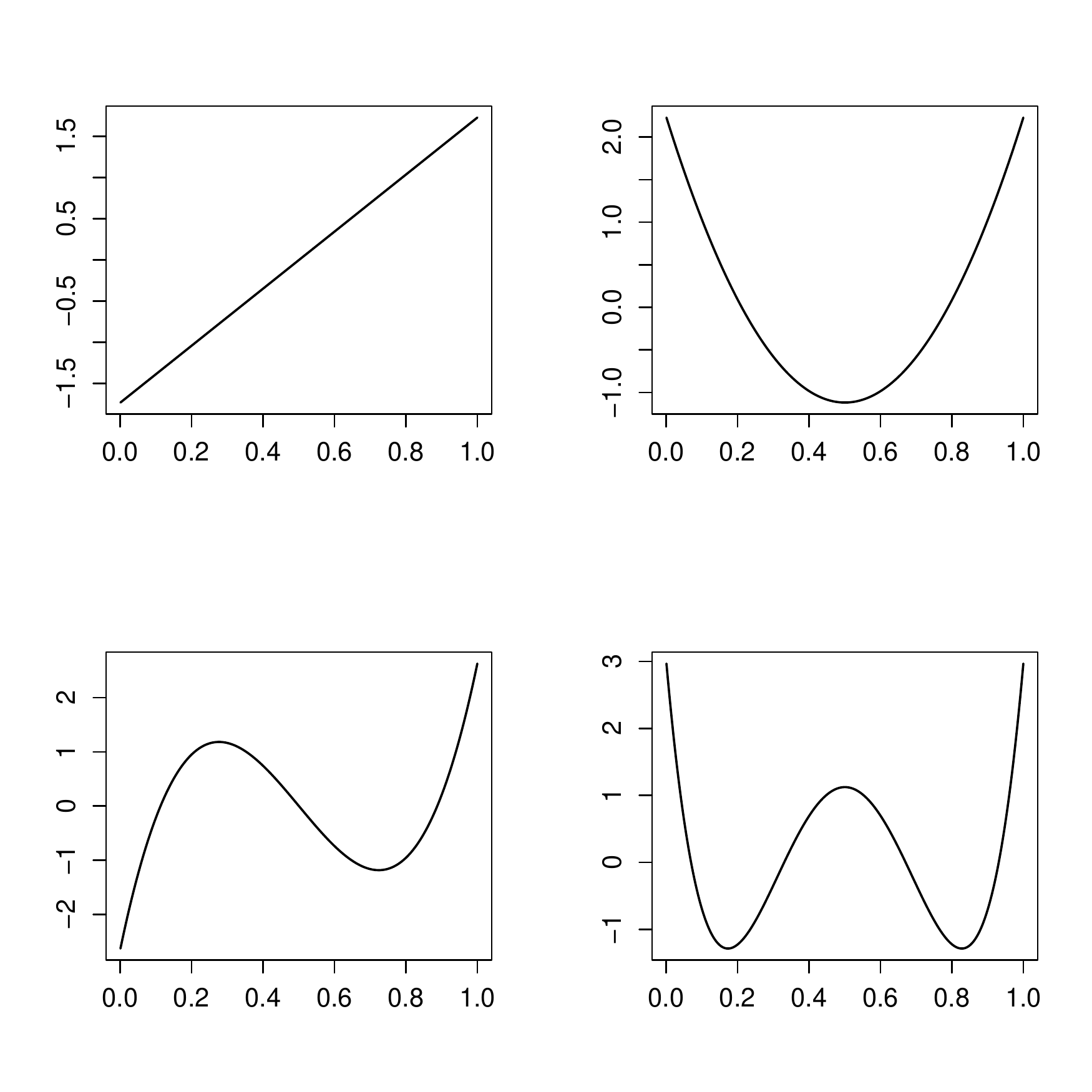}
\vspace{-.1em}
\caption{Shapes of first four orthonormal shifted Legendre polynomials on $[0,1]$. } \label{fig:leg}
\end{figure}

\subsection{Relations with Legendre polynomials}

For $X$ continuous, LP unit score functions yield Legendre polynomials, as shown in Fig \ref{fig:leg}. For $X$ discrete-uniform or for a sample distribution with all distinct values, LP score functions take the form of discrete Legendre polynomials, shown in Fig \ref{fig:score}(b). General piecewise constant shape of orthonormal LP unit score functions for $X$ discrete random variable is shown in Fig \ref{fig:score}(a), based on $n=250$ simulated samples from poisson with $\la=2$. Note that sample distributions are always discrete.

\section{LP Moments, LP Comoments}
The L moments method of data analysis, introduced by \cite{hos90}, was pioneered by \cite{sillitto1969} in the paper ``Derivation of approximants to the inverse distribution of a continuous univariate population from the order statistics of a sample.'' The goal is nonparametrically estimate distribution function $F(x;X)$ of variable $X$ in a way that is systematic (not ad hoc) and avoids high powers of $X$ (moments). \cite{sillitto1969} shows that continuous $X$ can be expressed as a polynomial of $F(X;X)$, the rank transform, or probability integral transform of $X$.  We extend this concept to \emph{general} $X$, discrete or continuous, by introducing the concept of LP moments and the multivariate generalization LP-comoments. We provide a compact representation and computationally efficient algorithm to compute the exact (unbiased version) sample estimates, which is known to be a challenging problem \citep{HosL}.

\subsection{LP Moments, Representation of X}
Define LP moments as $\LP(j;X)=\Ex[X T_j(X;X)]$. An equivalent representation using LP unit score functions (by expressing in quantile domain) is
\beq
\LP(j;X)=\int_0^1 Q(u;X) S_j(u;X) \dd u,
\eeq
which leads to the following important representation theorem by Hilbert space theory that provides nonparametric identification of models for \emph{mixed} $Q(u;X)$ by estimating its LP moments.
\begin{thm} \label{thm:qex}
Quantile function $Q(u;X)$ for a random variable $X$ with finite variance can be expressed as an infinite linear combination of LP unit score functions with coefficients as LP-moments
\beq \label{eq:qex}
Q(u;X) - \Ex(X) = \sum_{j>0} \LP(j;X)\, S_j(u;X),
\eeq
\end{thm}
From Theorem \ref{thm:qex} we derive the the following LP representation result of $X$ as a function of $F(X;X)$; It is similar to a Karhunen-Lo{\'e}ve expansion leading to a sparse representation of $X$ in a custom built orthonormal basis.
\begin{thm} \label{thm:Xlp}
A random variable $X$ (discrete or continuous) with finite variance admits the following decomposition with probability $1$
\beq  \label{eq:Xlp}
X - \Ex(X) = \sum_{j>0} \LP[j;X]\, T_j(X;X),
\eeq
\end{thm}
We call \eqref{eq:Xlp} the LP representation of $X$.
Theorem \ref{thm:Xlp} is a direct consequence of applying the following fundamental theorem to (\ref{eq:qex}).
\begin{lemma} \label{lemma:qf1}
For a random variable $X$ with  probability $1$,  $Q(F(X;X);X)=X$.
\end{lemma}
A published proof of this known fact is difficult to find; we outline proof in Appendix A2.
Theorem \ref{thm:Xlp} leads to a very useful fact about variance of $X$ given in the following corollary.
\begin{coro}
Variance decomposition formula in terms of LP moments
\[\var(X) = \sum_{j>0}\big|\LP(j;X)\big|^2\]
\end{coro}

\begin{rem}[Interpretation of LP shape statistic]
For a sample of $X$ continuous with distinct values one can derive the following representation of $\LP[1;\tF]$ as a linear combination of order statistics
\[\LP[1;\tF]\,=\, \dfrac{2 \sqrt{3}}{n\sqrt{n^2-1}} \sum_{i=1}^n\, X_{(i)} \Big( i\, - \,\dfrac{1}{2}(n+1) \Big).  \]
Our $\LP(1;X)$ is a measure of scale that is related to the Downton's estimator \citep{downton1966} and can be interpreted as a modified \emph{Gini mean difference coefficient} \citep{david1968}. Measures of skewness and kurtosis are $\LP(2;X)$ and $\LP(3;X)$.
\end{rem}
Table 2 in Appendix A3 reports numerical values of LP moments for some standard discrete and continuous distributions. For the Ripley data (described in the next section), the first four sample LP moments of $X=$ Age are $ 4.74, 1.49, 0.27, 0.11$ and $\Var(X)=24.84$. For $Y=$ GAG level, the first four sample LP moments are $8.08, 2.90, 1.89,1.05$ and $\Var(Y)=80.87$. The nonparametric robust LP moments provide a quick understanding of the shape of the marginal distributions, as depicted in Fig \ref{fig:marginalgag}.

\vskip.4em
An application of LP moment to test of normality is given by the following result. For $X \sim \cN(\mu,\si^2)$ we have
\beq \LP(1;X)\,=\, \si\, \Ex[X \Leg_1(\Phi(X));X] \,=\, 2\sqrt{3}\, \si \,\Ex[ \phi(X); X]\,=\,\sqrt{3/\pi}\, \si.\eeq
This suggests that D'Agostino's normality test \citep{dagos1971} can be computed as
\beq \Cor(X,F(X)) ~=~\LP[1;X]/\si.   \eeq

Define $X$ to be short tailed if $|\LP(1;\cZ(X)]|^2 > .95$, where $\cZ(X)=(X-\Ex[X])/\si$.

\subsection{LP Comoments, Dependence Analysis}
To model dependence and joint distribution of $(X,Y)$, define LP comoments as the cross-covariance between higher-order orthonormal LP score functions $T_j(X;X)$ and $T_k(Y;Y)$,
\beq \LP[j,k;X,Y]~=~\Ex[T_j(X;X)\,T_k(Y;Y)] ~~~\mbox{for}~ j,k>0.\eeq
From the LP-representation \eqref{eq:Xlp} of $X$ and $Y$ we obtain following expansion of covariance based on LP-comoments and LP-moments.
\begin{thm}
Covariance between $X$ and $Y$ admits the following decomposition
\beq
\Cov(X,Y)~=~\sum_{j,k>0}\LP(j;X) \LP(k;Y) \LP(j,k;X,Y)
\eeq
\end{thm}
Table 2 in Appendix A4 lists LP comoments of bivariate standard normal with correlation $\rho$ equal $0,.5,.9$.

\subsection{Ripley GAG Urine Data, Fisher Hair and Eye Color Data}

We will be using the following two interesting data sets to demonstrate how our LP algorithms answer scientific questions.

\begin{exmp}[Ripley GAG Urine Data]
This data set was introduced by \cite{ripley04}. The study collected the concentration of a chemical GAG in the urine of $n=314$ children aged from zero to seventeen years. The aim of the study was to produce a chart to help paediatricians to assess if a child's GAG concentration is `normal.' In Section 7.4 we present our comprehensive analysis of this data set.
Eq. \eqref{eq:glpco} computes the empirical LP-comoment matrix of order $m=4$ between $X$ Age and $Y$ GAG level for the Ripley data:
\beq \label{eq:glpco}
\widehat{\LP}\big[X ,Y \big]~=~\begin{bmatrix}
                              -0.908^{*}  &-0.010  &0.011 &0.035 \\[.38em]
                              0.032  &0.716^{*} &-0.071 &0.028\\[.38em]
                              0.064  &0.015 &-0.590^{*} &0.117\\[.38em]
                              -0.046 &-0.085 &-0.060 &0.425^{*}\\[.38em]
                             \end{bmatrix}\eeq

\end{exmp}

\begin{exmp}[Fisher Hair and Eye Color Data]
A two-way continency table describes the hair and eye color of $n=5387$ Scottish children, first analyzed by \cite{fisher1940}. Section 5.2 presents a comprehensive modeling of Fisher's data. For the discrete Fisher data (described in Section 5.2) the estimated LP-comoment matrix is given by
\beq \label{eq:glpco2}
\widehat{\LP}\big[X ,Y \big]~=~\begin{bmatrix}
                               0.423^{*} &0.024 &0.039 &-0.009 \\[.38em]
                              0.115^{*} &0.157^{*} &0.001 &-0.021\\[.38em]
                               -0.050 &0.085 &0.017 &-0.032\\[.38em]
                             \end{bmatrix}\eeq
\end{exmp}

\subsection{LP[1,1] as Generalized Spearman Correlation}
How to define a properly normalized Spearman correlation for discrete data with ties is an ongoing question. Recently, there has been various proposals discussed in \cite{denuit05,nevslehova07,mesfioui2010,blumentritt2012}, and \cite{genest2013}. It has been recognized that the range of Spearman for discrete data with ties depends on marginals and typically does not reach the bounds $\pm 1$; a surprising example is given in \cite{genest2007}.

\begin{exmp}[The Genest Example]
$(X,Y)$ Bernoulli random variables with marginal and joint distribution $\Pr(X=0)=\Pr(Y=0)=\Pr(X=0,Y=0)=p \in (0,1)$, which implies $Y=X$ almost surely, still the traditional $R_{{\rm Spearman}}(X,Y)=p(1-p)<1$. Similarly, if $\Pr(X=0)=p=1-q=\Pr(Y=1)$ and $\Pr(X=0,Y=1)=p \in (0,1)$, then $Y=1-X$ almost surely, yet $R_{{\rm Spearman}}(X,Y)=-p(1-p)>-1$.
\end{exmp}

Our $\LP[1,1;X,Y]$ (linear-rank) statistic not only resolves these inconsistencies (built-in adjustment based on mid-distribution based score polynomials $T_1[\Fm(X;X)]$ and $T_1[\Fm(Y;Y)]$) but also provides \emph{one single computing formula} for Spearman correlation that is applicable to discrete and continuous marginals (\emph{without} any additional adjustments). A straightforward calculation shows that $\big|\widehat{\LP}(1,1;X,Y)\big|=1$ for the Genest Example.
\vskip.2em

Moreover, our LP-comoment-based proposal \emph{systematically produces other higher-order nonlinear (ties-corrected) dependence measures}, therefore providing much more precise information about the dependence structure. A detailed investigation of this idea is beyond the scope of the current paper and will be discussed elsewhere.
\section{LP Skew Density Estimation}
For $X$ continuous, kernel density estimation \citep{parzen1962} is a popular nonparametric density estimation tool.
In this section we introduce a general nonparametric density estimation algorithm, which \emph{simultaneously} treats discrete and continuous data; we call it LP Skew density. We illustrate our method using Buffalo Snowfall data and Online rating discrete data from \texttt{Tripadvisor.com}. We further show how LP skew estimator can lead to a general algorithm for goodness-of-fit that unifies Pearson's chisquare \citep{pearson1900}  for $X$ discrete and Neyman's approach \citep{Neyman37} for $X$ continuous. Empirical application towards large-scale multiple testing is also discussed.
\subsection{Definition and Properties}
Probability density $f(x;X)$ of continuous $X$, or probability mass function of $p(x;X)$ of discrete $X$ are nonparametrically estimated by a model we call a LP Skew density, constructed as follows:
\vskip.25em
{\bf Step A} Choose suitable distribution $G(x)$, with quantile $Q(u;G)$, with probability density or probability mass function denoted $g(x)$. \vskip.25em
\vskip.25em
{\bf Step B} Define $d(u), 0<u<1$, by $d(G(x))=f(x;X)/g(x)$ or $d(G(x))=p(x;X)/g(x)$.

We call $d(u)$ comparison density; more explicitly, it is defined as
\beq d(u;G,F)=f(Q(u;G);X)/g(Q(u;G)), ~~ {\rm or}~~
d(u;G,F)=p(Q(u;G);X)/g(Q(u;G)).
\eeq
\vskip.25em
{\bf Step C}  An LP Skew density estimator of the true probability law is provided by an estimator
\beq
\widehat f(x)\,=\,g(x)\,\, \dhat(G(x)).
\eeq
\vskip.25em
An $L^2$ estimator of $d(G(x))$, or maximum entropy (exponential family) estimator of $\log d(G(x))$, is obtained by representing it as a linear combination of $T_j(x;G)$; LP score functions of G can be considered ``sufficient statistics'' for a parametric model of $d(u)$.
To select significant score functions, apply Akaike AIC type model selection theory or data driven density estimation theory of \cite{kall99} to select significant ``goodness of fit components,'' a fundamental concept defined as
\beq
\LP(j;G,F)~=~ \Ex[T_j(X;G);F]~=~\int_0^1 S_j(u;G) \dd D(u;G,F) ~=~\big \langle d, S_j  \big\rangle.
\eeq
\begin{thm}
Empirical process theory shows that when $F$ equals sample distribution of sample of size $n$ and $G$ equals true distribution $F$, then $\sqrt{n} \LP(j;G,F)$ is asymptotically $\cN(0,1)$ and provides a test of goodness of fit of $G$ to true $F$.
\end{thm}

\subsection{Applications}
\begin{figure*}[!thb]
 \centering
 \includegraphics[height=.38\textheight,width=.47\textwidth,keepaspectratio,trim=.5cm .5cm .5cm .5cm]{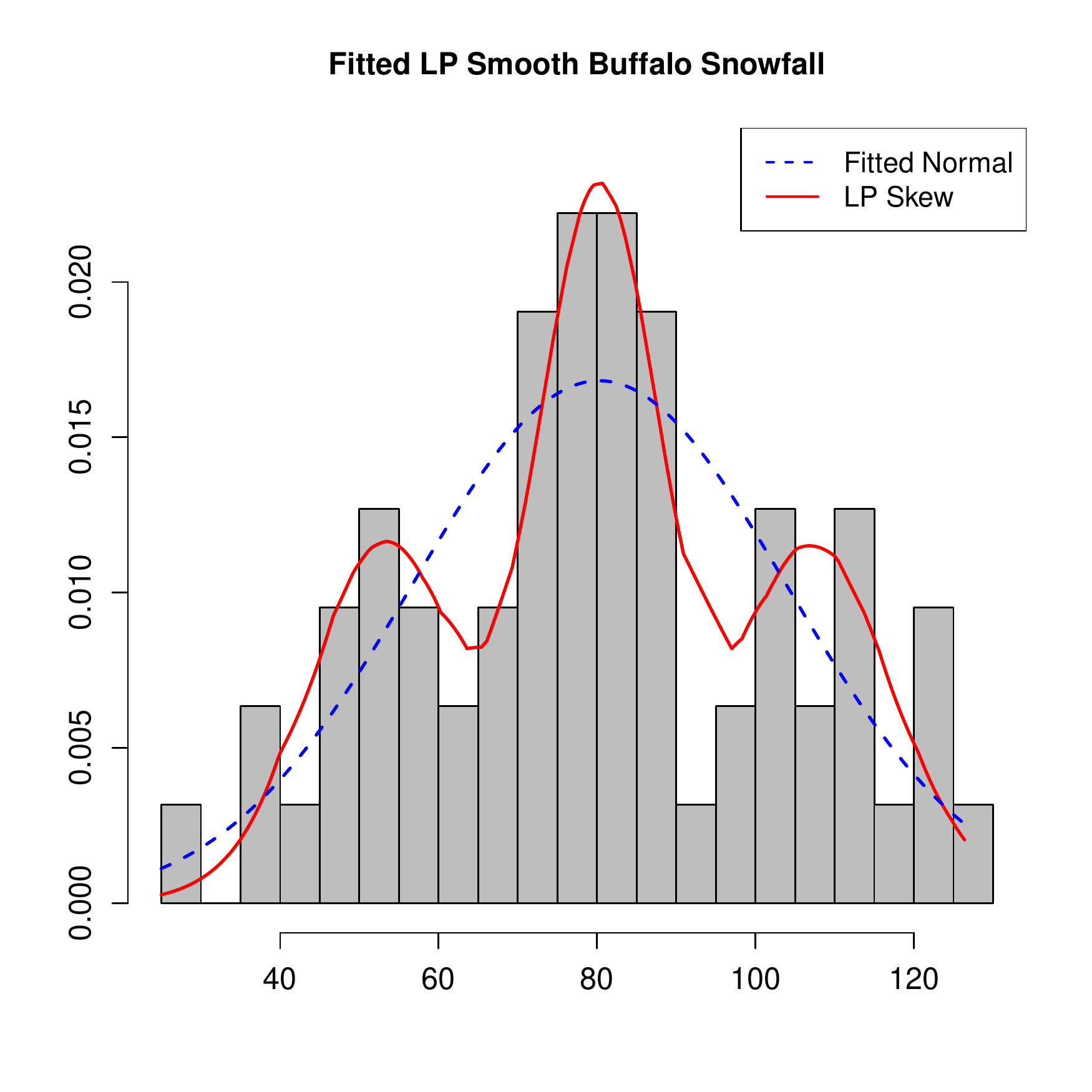}
 \includegraphics[height=.38\textheight,width=.47\textwidth,keepaspectratio,trim=.5cm .5cm .5cm .5cm]{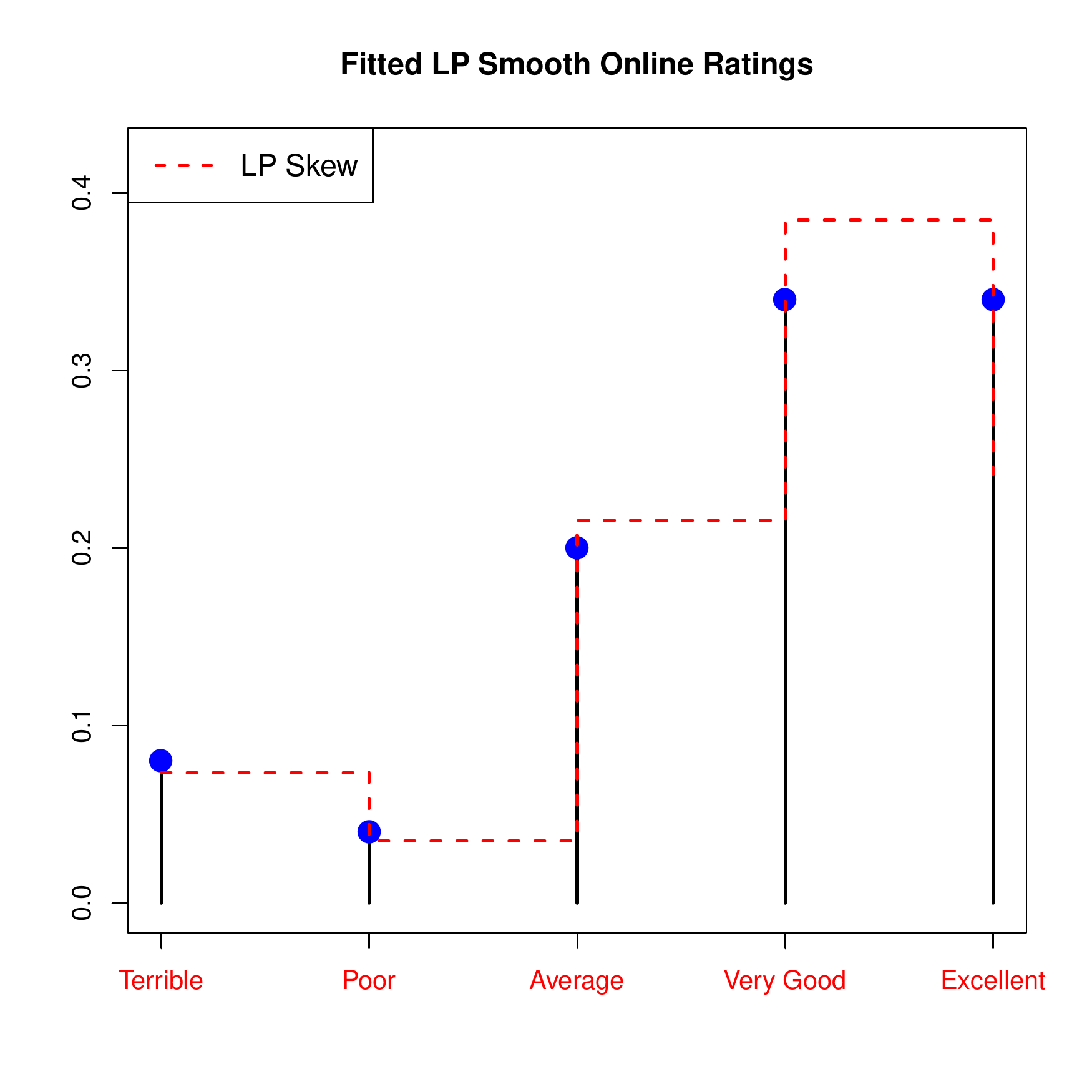}
\vspace{-.1em}
\caption{LP Skew density estimates. (a) Left; Buffalo Snowfall data. The baseline measure $\cN(\mu=80.29,\si=23.72)$ shown in blue dotted line; (b) Right; Online rating data with base line measure as Poisson distribution with $\la=3.82$. } \label{fig:LSD}
\end{figure*}

\subsubsection{Smooth Density Modeling}
Two examples are given to illustrate how LP skew smoothing methods nonparametrically identify the underlying probability law.
\vskip.3em

Figure \ref{fig:LSD}(a) shows the Buffalo Snowfall data, previously analyzed by \cite{parzen79}. We select the normal distribution (with estimated mean $80.29$ and standard deviation $23.72$) as our baseline density $g(x)$. Data-driven AIC selects $\dhat(u)=1-.337 S_6(u)$. The resulting LP Skew estimate shown in Fig \ref{fig:LSD}(a) strongly suggests tri-modal density.
\vskip.3em

Fig \ref{fig:LSD}(b) shows online rating discrete data from \texttt{Tripadvisor.com} of the Druk hotel with baseline model $g(x)$ being Poisson distribution with fitted $\la=3.82$. Data-driven AIC selects $\dhat(u)=1-.55 S_2(u)-.4 S_3(u)+.23 S_4(u)$. The estimated LP Skew probability estimates $\widehat p(x;x)$ captures the overall shape, including the two modes, which reflects the mixed customer reviews: (i) detects the sharp drop in probability of review being `poor' (elbow shape); (ii) both peaks at the boundary, and (iii) finally, the ``flat-tailed'' shape at the end of the two categories (`Very Good' and `Excellent'). Compare our findings with Fig 12 of \cite{Shmueli2013}, where they have used a substantially complex and computationally intensive model (Conway-Maxwell-Poisson mixture model), still fail to satisfactorily capture the shape of the discrete data.

%
%
%

\subsubsection{General Goodness-of-fit Principle}

LP Skew modeling naturally leads to a framework for \emph{distribution-free} goodness-of-fit procedure, valid for \emph{mixed} $X$ (discrete or continuous). A distribution-free goodness-of-fit framework for discrete data is comparatively much less developed and more challenging than continuous data \citep{d1986}. Some recent attempts have been made by \cite{best1999,best2003} and \cite{khmaladze2013}.

To test whether the true underlying distribution $F$ equals $G$, we follow the following procedure:
\vskip.35em

{\bf Step A} Estimate the $\dhat(G(x))\,=\,1+\sum_j \LP[j;G,\tF]\, T_j(x;G)$.
\vskip.25em

{\bf Step B} Compute the goodness-of-fit test statistic $\int_0^1 \widehat{d}^2\,-\,1 \,=\,\sum_j  \big| \LP[j;G,\tF]\big|^2$. Justify the form of the test statistic by recognizing that under null the comparison distribution limit process $\sqrt{n} \big[\wtD(u;G,F)-u\big] \xrightarrow{D} B(u)$, Brownian Bridge for $0<u<1,$ with RKHS norm squared $\|h\|^2=\int_0^1 |h'(u)|^2 \dd u$.
\vskip.35em

The following result shows one remarkable equality between LP goodness-of-fit statistic and Chi squared statistic for discrete data.

\begin{thm}
For $X$ discrete with probable values $x_j$ and observe sample size $n$, one compares sample probabilities $\tilde p(x)$ with population model $p_0(x)$ by Chi-Squared statistic, which we represent by $n \CHIDIV$. Then LP goodness of fit statistics numerically equals to $\CHIDIV\big[p_0 \,\| \,\tilde p  \big]$.
\end{thm}
To prove the result, verify the following:
\beq \mbox{CHIDIV}\big[\,p_0 \,\| \,\tilde p  \big]\, = \,\sum_x p_0(x)\big[\tilde p(x)/p_0(x) \,-\,1\big]^2\,=\,\int_0^1 \big[d(u;p_0;\tilde p)-1\big]^2 \dd u.\eeq

Chi-square statistic is a ``raw-nonparametric'' information measure,  which we interpret by finding an approximately equal (numerically) ``smooth" LP information measure with \emph{far fewer degrees of freedom} because it is the sum of squares of only a few data-driven $\LP[j;G,\tF]$. LP goodness-of-fit statistic for discrete distribution can be interpreted as Chi-square statistic with the \emph{fewest possible degrees of freedom} that makes it more powerful. Our approach utilizes the automatic and universal algorithm to custom-construct orthonormal score polynomials of $F(X;X)$ for each given $X$; traditional methods construct orthonormal polynomials on a \emph{case-by-case basis} by each time solving the heavy-duty Emerson recurrence \citep{emerson1968,best1999,best2003} relation, which could often be quite complicated for non-standard distributions.

%

\begin{figure*}[!thb]
 \centering
 \includegraphics[height=.38\textheight,width=.47\textwidth,,keepaspectratio,trim=1cm 1cm .5cm .5cm]{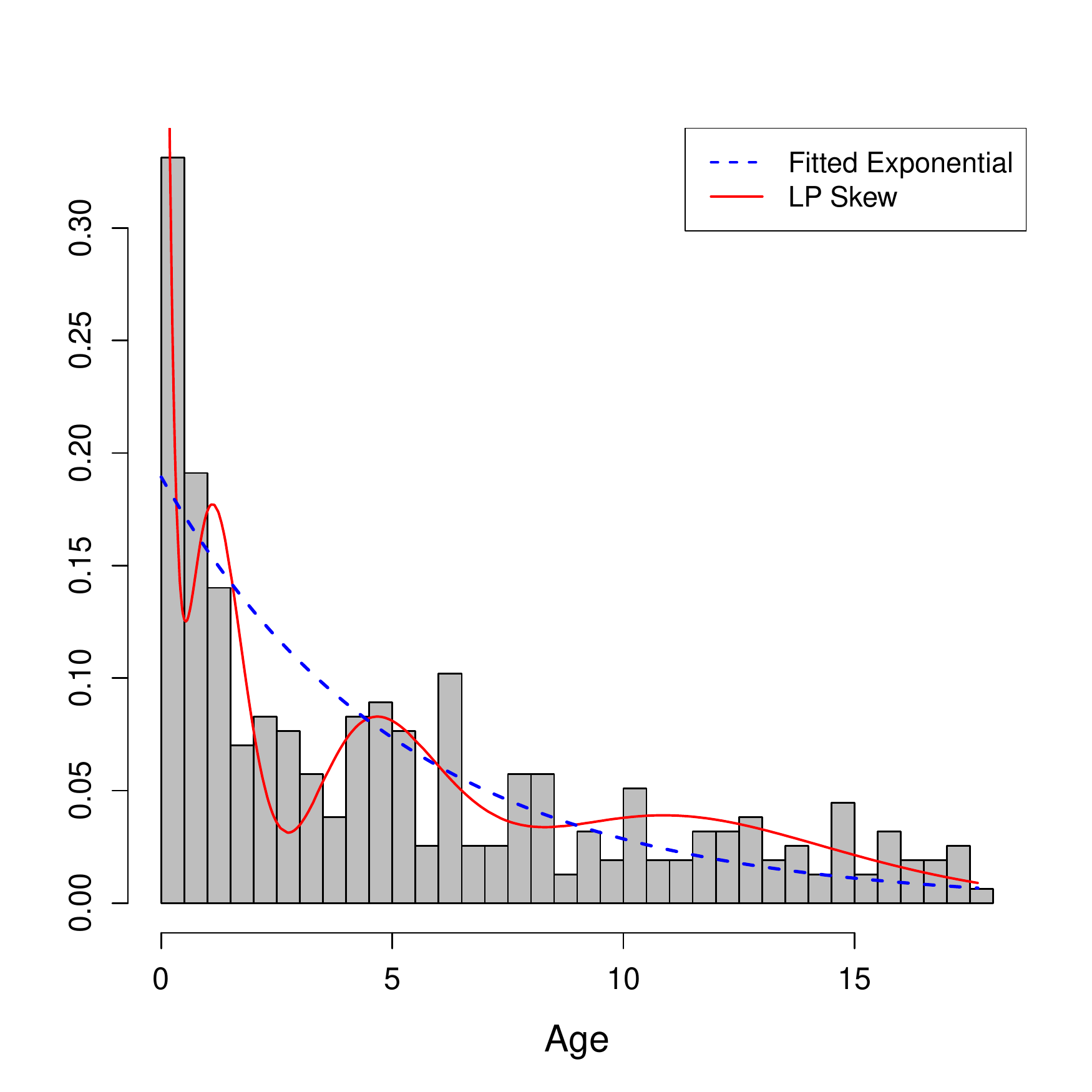}
 \includegraphics[height=.38\textheight,width=.47\textwidth,,keepaspectratio,trim=.5cm 1cm .5cm 1cm]{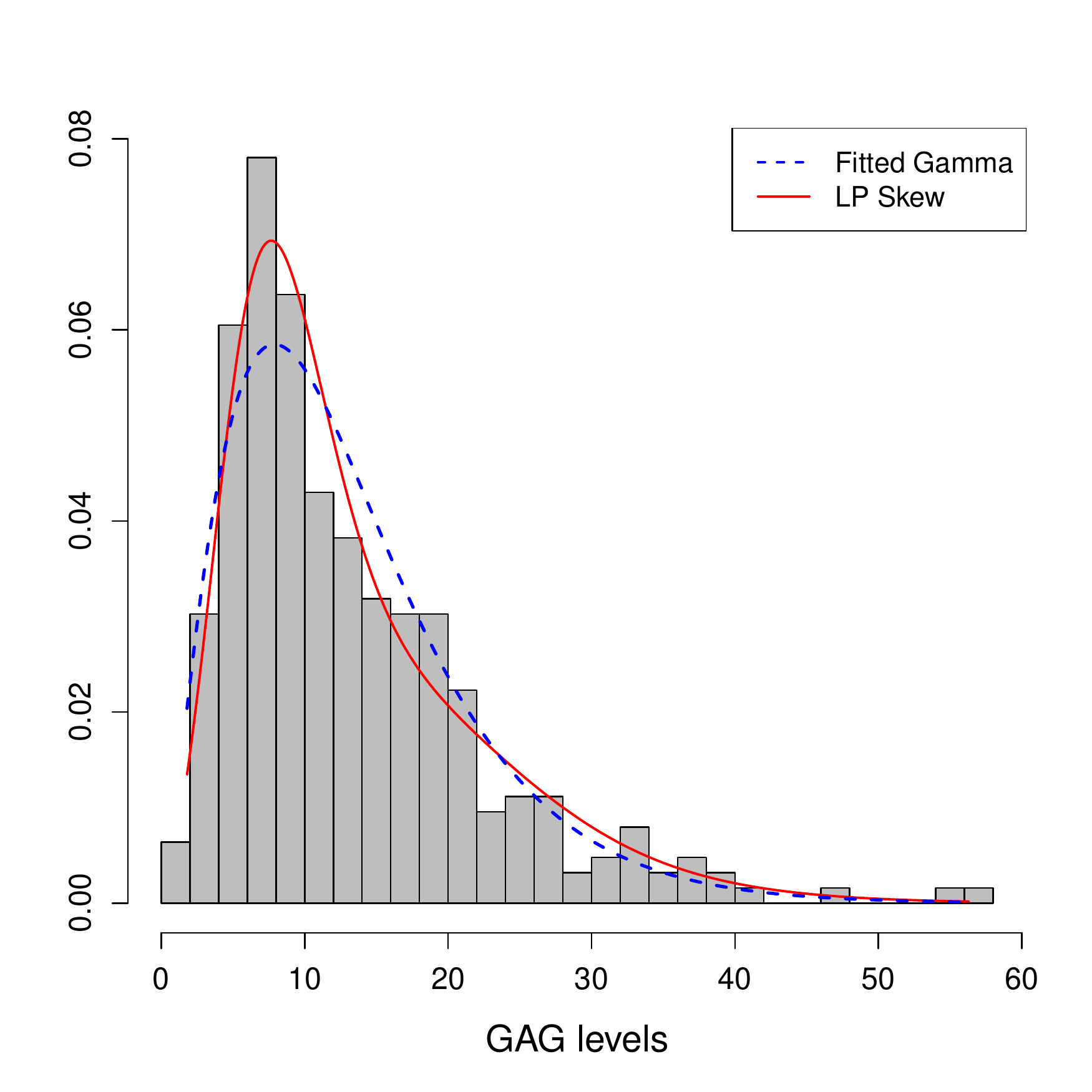}
\vspace{-.1em}
\caption{Goodness-of-fit of the marginal distributions of Ripley data.} \label{fig:marginalgag}
\end{figure*}
\vskip.2em
{\bf Examples.} For $Y=$ GAG level in Ripley data, we would like to investigate whether the parametric model Gamma fits the data. We select the Gamma distribution with shape $2.5$ and rate parameter $.19$ estimated from the data as our base line density $g(x)$. Estimated smooth comparison density is given by $\dhat(u;G,F)=1+.16 S_3(u)$, which indicates that the underlying probability model is slightly more skewed than the proposed parametric model. Fig \ref{fig:marginalgag}(b) shows the final LP Skew estimated model.

For $X=$ AGE in Ripley data, we choose exponential distribution with estimated mean parameter $ 5.28$ as our $g(x)$. The estimated  $\dhat(u;G,F)=1+.32S_2(u)-.25S_3(u)-.22S_5(u) -.39S_7(u)$ with LP goodness-of-fit statistic as $\int \dhat^2 -1= 0.365422$. Fig \ref{fig:marginalgag}(a) shows that the exponential model is inadequate for capturing the plausible multi-modality and heavy-tail of age distribution.

\subsubsection{Tail Alternative \& Large-Scale Multiple Testing} \label{sec:ta}
\begin{figure*}[!thb]
 \centering
 \includegraphics[height=.5\textheight,width=.80\textwidth,trim=1.5cm .5cm 1.5cm .5cm]{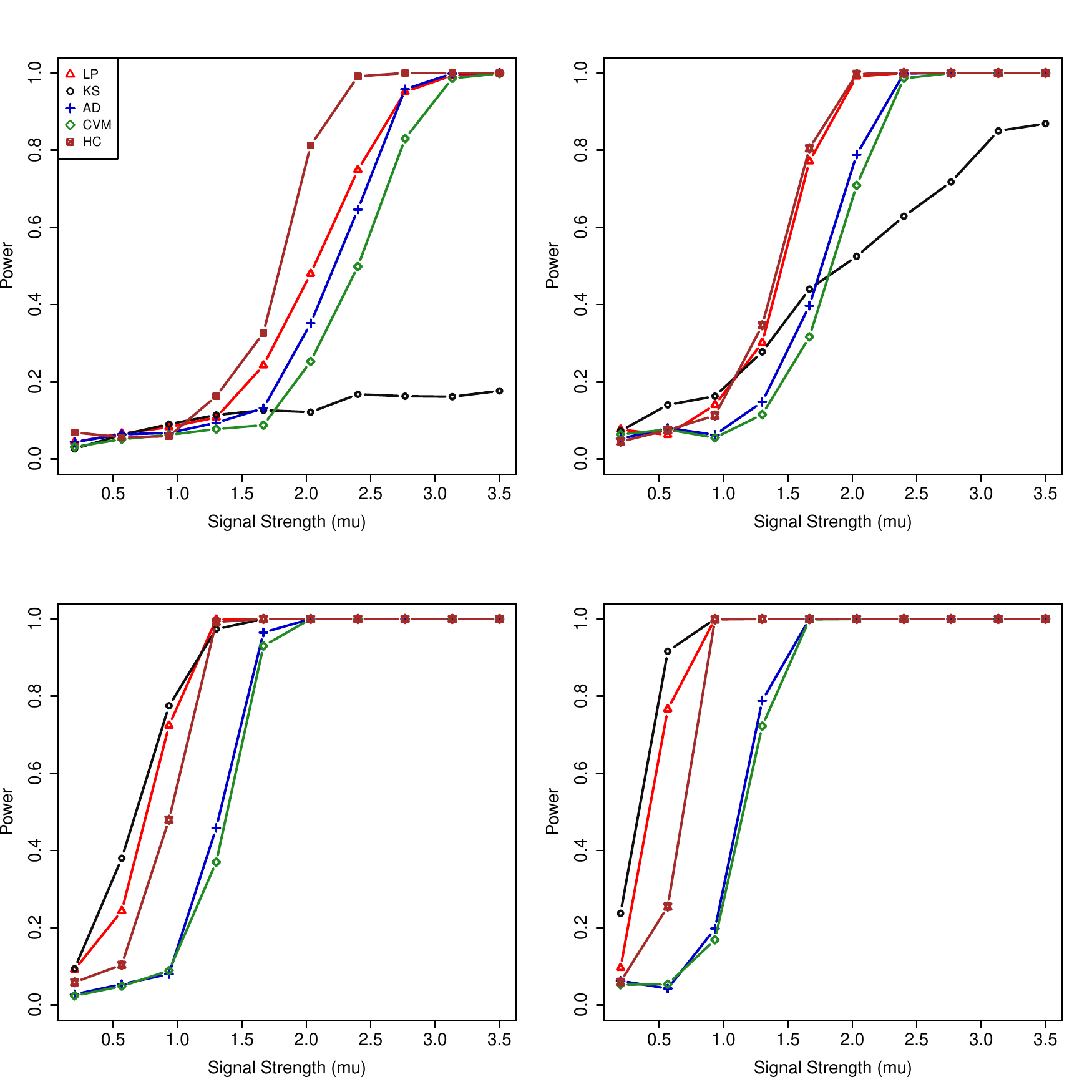}
\vspace{-.1em}
\caption{Power comparison of LP goodness-of-fit statistic under tail alternative (Section \ref{sec:ta}). The sparsity level changes $\pi =\{.01,.02,.05,.09\}$. } \label{fig:GOF}
\end{figure*}

We investigate the performance of the LP goodness-of-fit statistic under `tail alternative' (detecting lack of fit ``in the tail''). Consider the testing problem $H_0:F=\Phi$ verses Gaussian contamination model $H_1:F=\pi \Phi + (1-\pi) \Phi(\cdot - \mu)$, for $\pi \in (0,1)$ and $\mu>0$, which is critical component of multiple testing. It has been noted that traditional tests like Anderson–Darling (AD), Kolmogorov–Smirnov (KS) and Cram{\'e}r-von Mises (CVM) lack power in this setting.  We compare the performance of our statistic with four other statistics: AD, CVM, KS and Higher-criticism \citep{HC08}. Fig \ref{fig:GOF} shows the power curve under different sparsity level $\pi=\{.01,.02,.05,.09\}$ where the signal strength $\mu$ varies from $0$ to $3.5$ under each setting, based on $n=5000$ sample size. We used 500 null simulated data sets to estimate the 95\% rejection cutoff for all methods at the significance level $0.05$ and used $500$ simulated data sets from alternative to approximate the power.
 LP goodness-of-fit statistic does a satisfactory job of detecting deviations at the tails under a broad range of $(\mu,\pi)$ values, which is desirable as the true sparsity level and signal strength are unknown quantities. Application of this idea recently extended to modeling local false discovery rate in \cite{D11d}.

\section{LP Copula Density Estimation}
Copula density, pioneered by \cite{sklar1959} for continuous marginals, provides a general tool for dependence modeling. The copula density estimation problem is known to be a highly challenging problem, especially for discrete data as noted in \cite{genest2007} ``A Primer on Copulas for Count Data.''
Copula Distribution at $u=F(x;X),v=F(y;Y)$ for some $x,y$ is defined as
\beq \Cop(u,v;X,Y)=F\big(Q(u;X),Q(v;Y);X,Y\big).\eeq
For discrete and continuous $X$ and $Y$ define copula density for $0<u,v<1$ as
\beq
\cop(u,v;X,Y)\,=\,d(v;Y,Y|X=Q(u;X))\,=\,d(u;X,X|Y=Q(v;Y)),~0<u,v<1.
\eeq
This section provides \emph{two fundamental nonparametric representations} of copula density that hold for discrete and continuous marginals at the same time.

\subsection{Representation and Estimation}
\begin{thm}[LP Representation of Copula Density] \label{thm:lpc}
Square integrable copula density admits the following representation as infinite series of LP product score functions
\beq \label{thm:lpcop}
\cop(u,v;X,Y)-1= \sum_{j,k>0}  \LP(j,k;X,Y)\, S_j(u;X)\, S_k(v;Y),~ 0<u,v<1.\eeq
Equivalently,
$$\int_{[0,1]^2} \dd u \dd v\, d\big[v;Y,Y|X=Q(u;X)\big] \,S_j(u;X)\, S_k(v;Y)\,=\,\LP(j,k;X,Y).$$
\end{thm}

A proof of copula LP representation is provided by representations of conditional copula density and conditional expectations:
\bea d(v;Y,Y|X=Q(u;X)) &=& \sum_k   S_k(v;Y)\, \Ex[T_k(Y;Y)|X=Q(u;X)] \\
\Ex[T_k(Y;Y)|X=Q(u;X)]&=& \sum_j   S_j(u;X)\, \Ex[T_j(X;X) T_k(Y;Y)]. \eea

The LP copula representation provides data driven `smooth' $L^2$ estimators of the copula density after constructing custom-built score functions and LP comoments.

\vskip.45em

We build the maximum entropy LP exponential copula density model based on the selected product score functions
\beq
\log \cop(u,v;X,Y)\,= \,\sum_{j,k>0} \te_{j,k}\,S_j(u;X)\,S_k(v;Y) \,-\,K(\te),\,\,~0<u,v<1.
\eeq

We recommend plotting the slices of the copula along with the copula density for better insight into the nature of the dependence. Plot conditional comparison densities $d(v;Y,Y|X=Q(u;X))$ as a function of $v$ for selected values of $u$, also $d(u;X,X|Y=Q(v;Y))$ as a function of $u$ for selected values of $v$. The slices describes the conditional dependence of $X$ and $Y$, which will be utilized for correspondence analysis (in Section 5) and quantile regression (in Section 7).

\begin{figure*}[!thb]
 \centering
 \includegraphics[height=\textheight,width=.47\textwidth,keepaspectratio,trim=2.5cm 2cm 1.5cm 2.5cm]{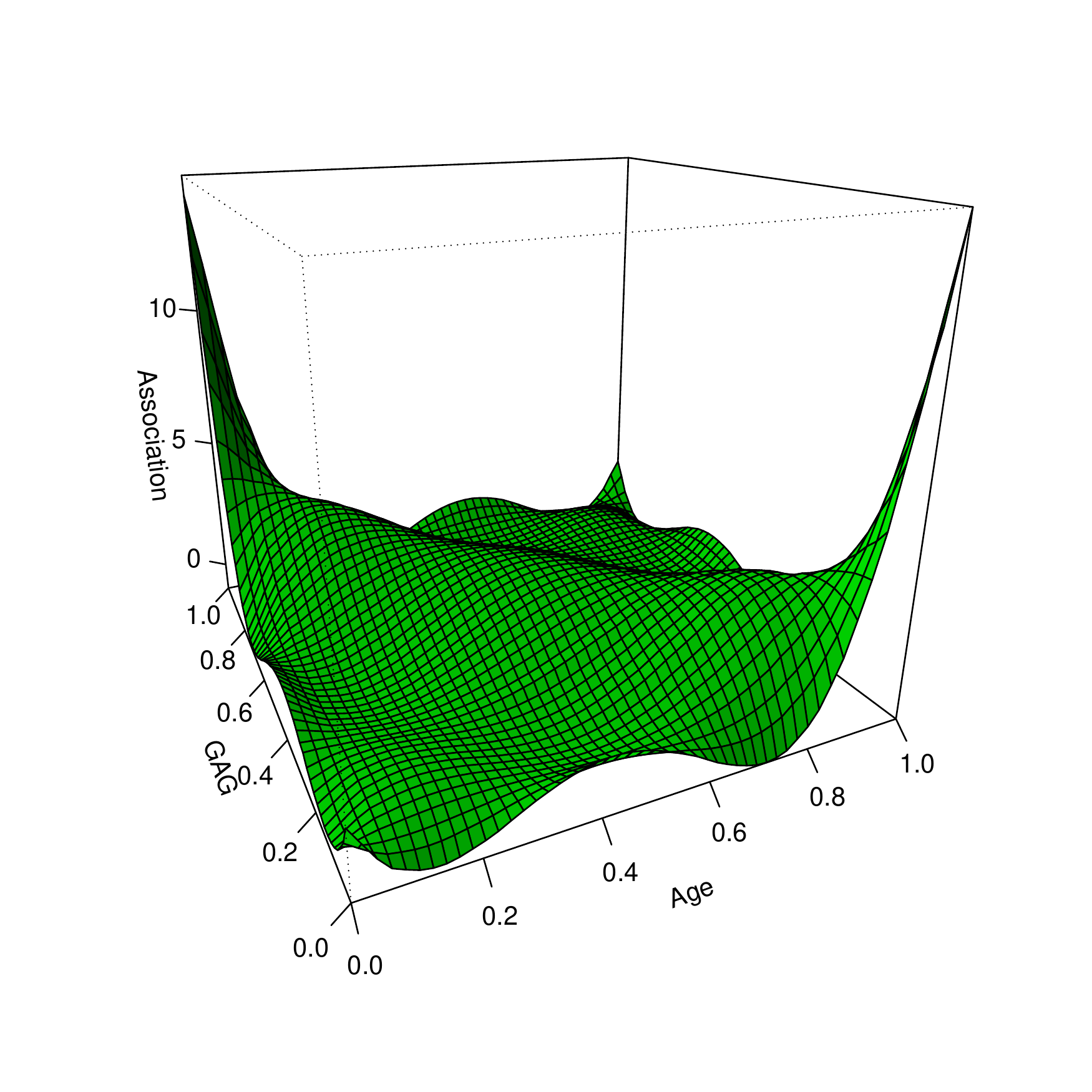}
  \includegraphics[height=\textheight,width=.49\textwidth,keepaspectratio,trim=1cm 2cm 2.5cm 2.5cm]{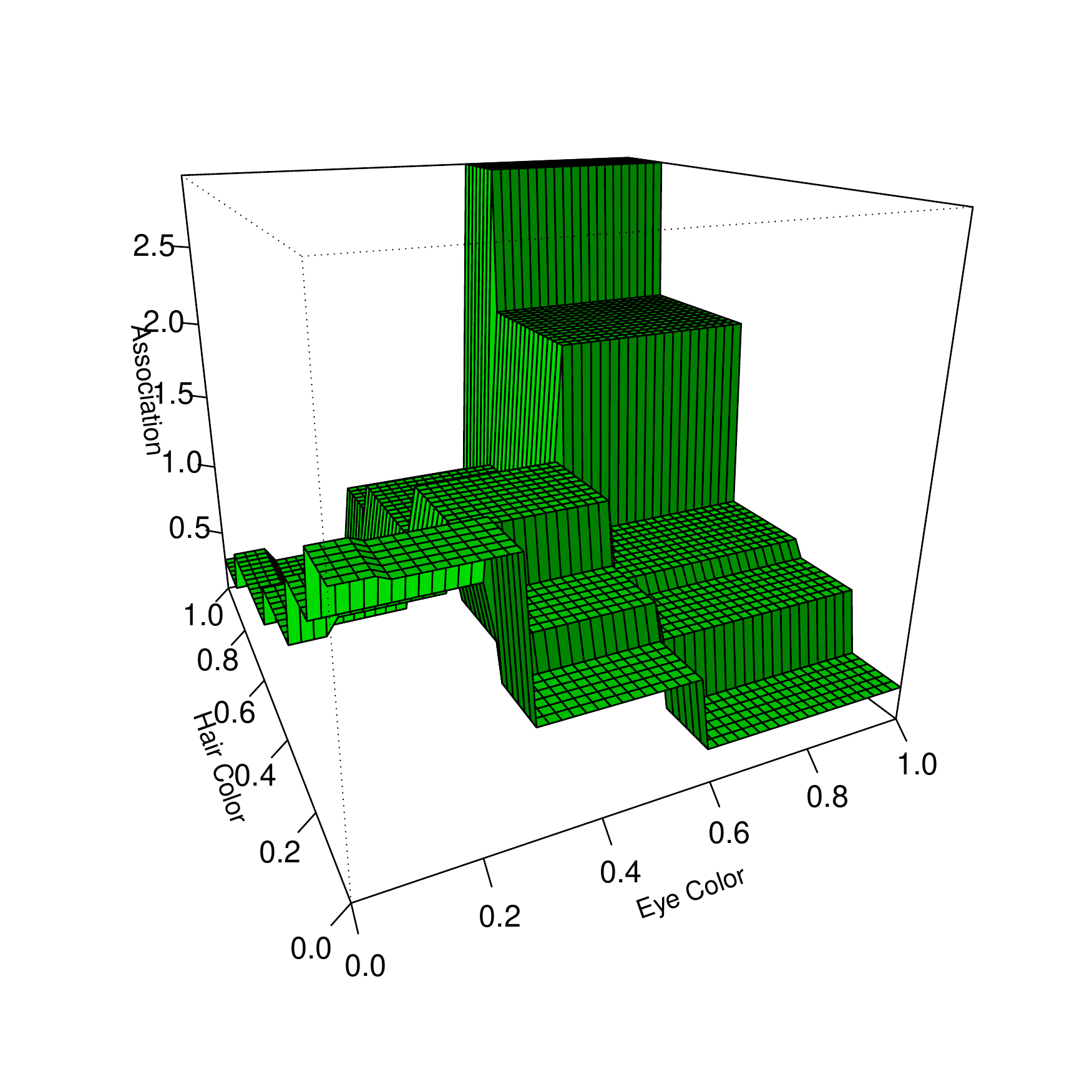}
\vspace{-.1em}
\caption{LP nonparametric copula density estimates for Ripley and Fisher data.} \label{fig:cop}
\end{figure*}

\subsection{Example}
Fig \ref{fig:cop}(a) displays the LP nonparametric copula estimate for the $(X,Y)$ Ripley data given by
\beq \widehat{\cop}(u,v;X,Y)\,=\,1 -0.91\,S_1(u)S_1(u)S_1(v)   + 0.71\,S_2(u)S_2(v)  - 0.59\,S_3(u)S_3(v)  +  0.42\,S_4(u)S_4(v)  \eeq
\vskip.35em

For $(X,Y)$ discrete Fisher data ($4 \times 5$ contingency table), the resulting smooth nonparametric checkerboard copula estimate is given by
\beq \widehat{\cop}(u,v;X,Y)\,=\,1 -0.423\,S_1(u)S_1(u)S_1(v)   + 0.157\,S_2(u)S_2(v) + 0.115\,S_2(u)S_1(v), \eeq
shown in Fig \ref{fig:cop}(b). Note the piece-wise constant bivariate discrete smooth structure of the estimated copula density.

\section{LP Canonical Copula Density Estimation}
We present two more copula expansions (Karhunen-Lo{\'e}ve type expansion) that are valid for both discrete and continuous marginals. For $(X,Y)$ discrete, we show that the vast literature on dependence modeling of two-way contingency table \citep{goodman1991,good96} can be elegantly unified by this formulation.
\subsection{Representation and Estimation} \label{thm:canc}
\begin{thm}[LP Canonical Expansion of Copula Density]
Square integrable copula density admits the following two representations ($L^2$ and maximum entropy exponential model) based on the sequence of singular values of LP-comoment kernel $\la_1 \ge \la_2 \ge \cdots$ for $0<u,v<1$,
\bea \label{eq:ccdecom} \cop(u,v)&=&1+ \sum_{k=1}^\infty\la_k \,\phi_k(u;X)\,\psi_k(v;Y), \\
\log \cop(u,v)&=&\sum_{k=1}^\infty \ga_k \,\phi_k(u;X)\,\psi_k(v;Y) \,-\,K(\ga).
 \eea
\end{thm}
Replace LP-comoment matrix appearing in Theorem \ref{thm:lpc} by its singular value decomposition (SVD) to get the required copula representation \eqref{eq:ccdecom}. Note that the singular values of the LP-comment matrix can be interpreted as the canonical correlation between the nonlinearly transformed random vectors $TX =(T_1(X;X), \ldots, T_m(X;X))'$ and $TY =(T_1(Y;Y), \ldots, T_m(Y;Y))'$. For this reason, we call this LP canonical copula representation.

For a $I \times J$ contingency table, the maximum number of terms in our LP copula expansion could be $\min(I,J)-1$. For $(X,Y)$ discrete Fisher data, the estimated LP canonical copula is given by $\cop(u,v)=1+.45\,\phi_1(u) \psi_1(v) + .17\,\phi_2(u) \psi_2(v)$.

\begin{table}[th!]
\caption{Fisher Data: Displayed joint probabilities $\Pr(X=x,Y=y)$ ; marginal probabilities $\Pr(X=x)$, $\Pr(Y=y)$; and corresponding mid-distributions, as well as coefficients of conditional comparison densities $d(v;Y,Y|X=x)$ and $d(u;X,X|Y=y)$.}
\vskip1em
\renewcommand{\tabcolsep}{.18cm}
\renewcommand{\arraystretch}{1}
\centering 
\begin{tabular}{c ccccccccc} 
\hline\hline 
&\multicolumn{5}{c}{\textbf{Hair Color}}&& \\
\textbf{Eye Color}  &Fair &Red &Medium &Dark &Black&$p(x;X)$&$\Fm(x;X)$ &$\la_1\phi_1$&$\la_2\phi_2$\\ [0.5ex]
\hline 
Blue & .061 & .007 &  .045 & .020&  .001&.14 &.07 &-0.400 &-0.165 \\ 
Light & .128 &.021  &.108 &.035 &.001 &.29 &.28 & -0.441 &-0.088\\
Medium & .064 &.016  &.168 &.076 &.005& .33&.59 &0.034  &0.245\\
Dark & .018 &.009  &.075 &.126 &.016&.24&.88  &0.703 &-0.134\\ 
$p(y;Y)$ &.27 &.05&.40&.26&.02    &&\\
$\Fm(y;Y)$ &.135  &.295&.52&.85&.99   &&  \\
$\la_1\psi_1$ &-0.544  &-0.233 &-0.042  &0.589  &1.094 \\
$\la_2\psi_2$ &-0.174   &-0.048 &0.208 &-0.104 &-0.286\\[1ex]
\hline \hline
\end{tabular} \label{tab:fisher}

\end{table}

\subsection{Correspondence Analysis: LP Functional Statistical Algorithm}
Given a $I \times J$ two-way contingency table, we seek to graphically display the association among the different row and column categories via correspondence analysis.

\subsubsection{Algorithm} LP correspondence analysis algorithm can be described as follows:
\vskip.35em
{\bf Step A} Estimate LP canonical copula density \eqref{eq:ccdecom} by performing SVD on LP-comoment matrix.
\vskip.25em

{\bf Step B} Define the LP principal profile coordinates $ \mu_{ik}=\la_k \phi_k(u_i;X) $ and $\nu_{jk}=\la_k \psi_k(v_j;Y)$.
\vskip.25em

{\bf Step C} Jointly display the row and column profiles together in the same plot to get the LP correspondence map.
\vskip.45em

\subsubsection{Example} 
\begin{figure}[tbh!]
 \centering
 \includegraphics[height=.45\textheight,width=\textwidth,keepaspectratio,trim=1.5cm .5cm 1.5cm .5cm]{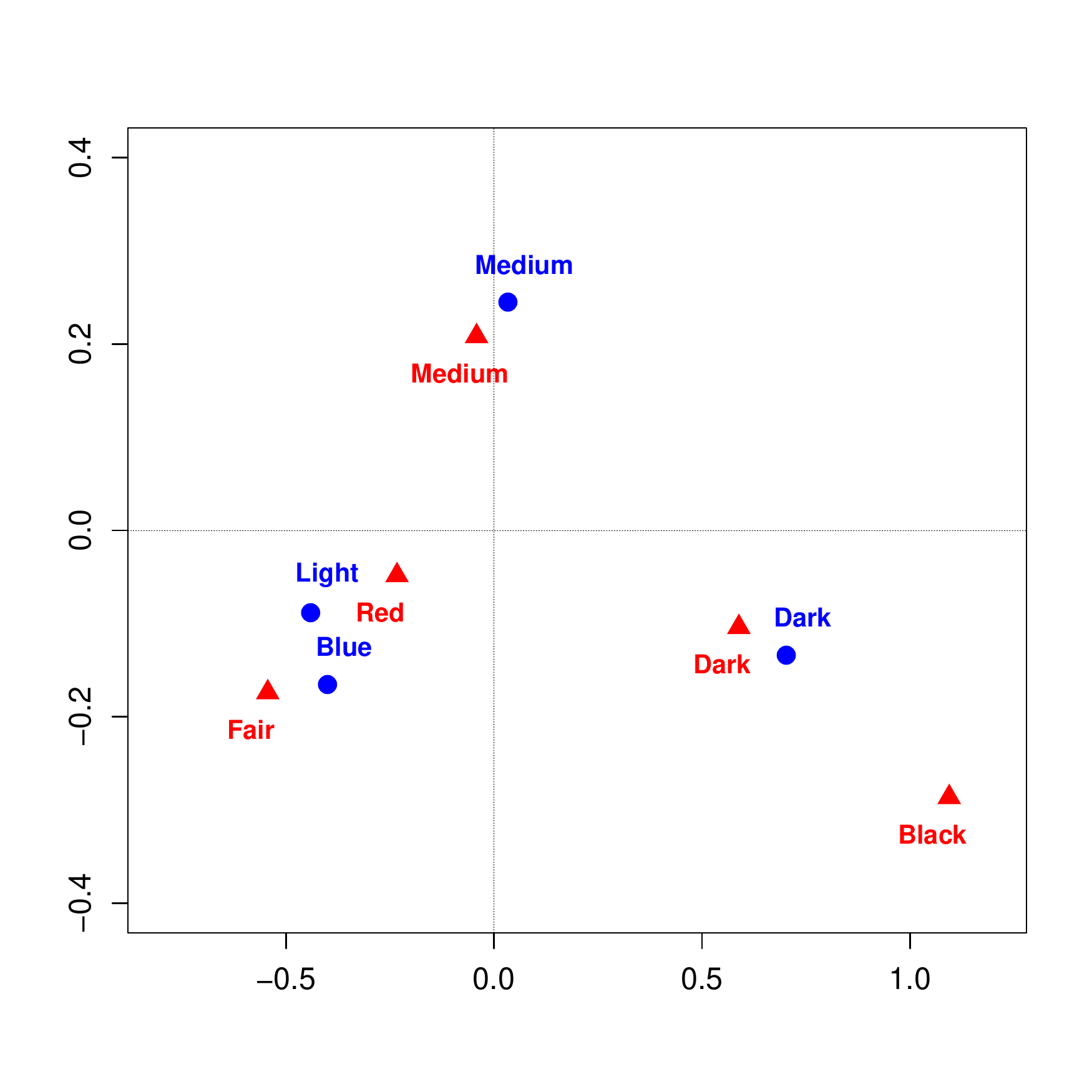}
\vspace{-.1em}
\caption{(color online) LP correspondence map for Fisher Hair and Eye color data}
\label{fig:fcorsp}
\end{figure}

Table \ref{tab:fisher} calculates the joint probabilities and the marginal probabilities of the Fisher data. The singular values of the estimated LP-comoment matrix are $0.446, 0.173, 0.029$. The rank-two approximation of LP canonical copula can adequately model (99.6\%) the dependence. Table \ref{tab:fisher} computes the row and column scores using the aforementioned algorithm. Fig \ref{fig:fcorsp} shows the bivariate LP correspondence analysis of Fisher data, which displays the association between hair color and eye color.
\vskip.35em

\subsubsection{LP Unification of Two Cultures} Our $L^2$ canonical copula-based row and column scoring (using LP algorithm) reproduces the classical correspondence analysis pioneered by \cite{hirschfeld1935} and \cite{benzecri1969}. Alternatively, we could have used the LP canonical exponential copula model to calculate the row and column profiles by $ \mu'_{ik}=\ga_k \phi_k(u_i;X) $ and $\nu'_{jk}=\ga_k \psi_k(v_j;Y)$, known as Goodman's profile coordinates pioneered in  \cite{goodman1981} and \cite{goodman1985}, which leads to the graphical display sometimes known as a weighted logratio map \citep{greenacre2010} or the spectral map (SM) \citep{lewi1998}. Our LP theory unifies the `two cultures' of correspondence analysis, Anglo multivariate analysis (usual normal theory) and French multivariate analysis (matrix  algebra data analysis), using LP representation theory of discrete checkerboard copula density.
\vskip.35em

Our $L^2$ canonical copula density (5.1) is sometimes known as correlation model or CA (Correspondence Analysis) model and the exponential copula model (5.2) as RC (Row-Col) association model \citep{goodman1991,good96}. Since our LP-scores satisfy: $\Ex[T_j(X;X)]=\sum_i T_j(x_i;X) p(x_i)=0$ and $\Var[T_j(X;X)]=\sum_i T^2_j(x_i;X) p(x_i)=1$ by construction, this association model also known as the \emph{weighted} RC model.

\vskip.35em

\subsubsection{Functional Statistical Interpretation}
\begin{figure}[tbh!]
 \centering
 \includegraphics[height=\textheight,width=\textwidth,keepaspectratio,trim=1.5cm .5cm 0cm 1.5cm]{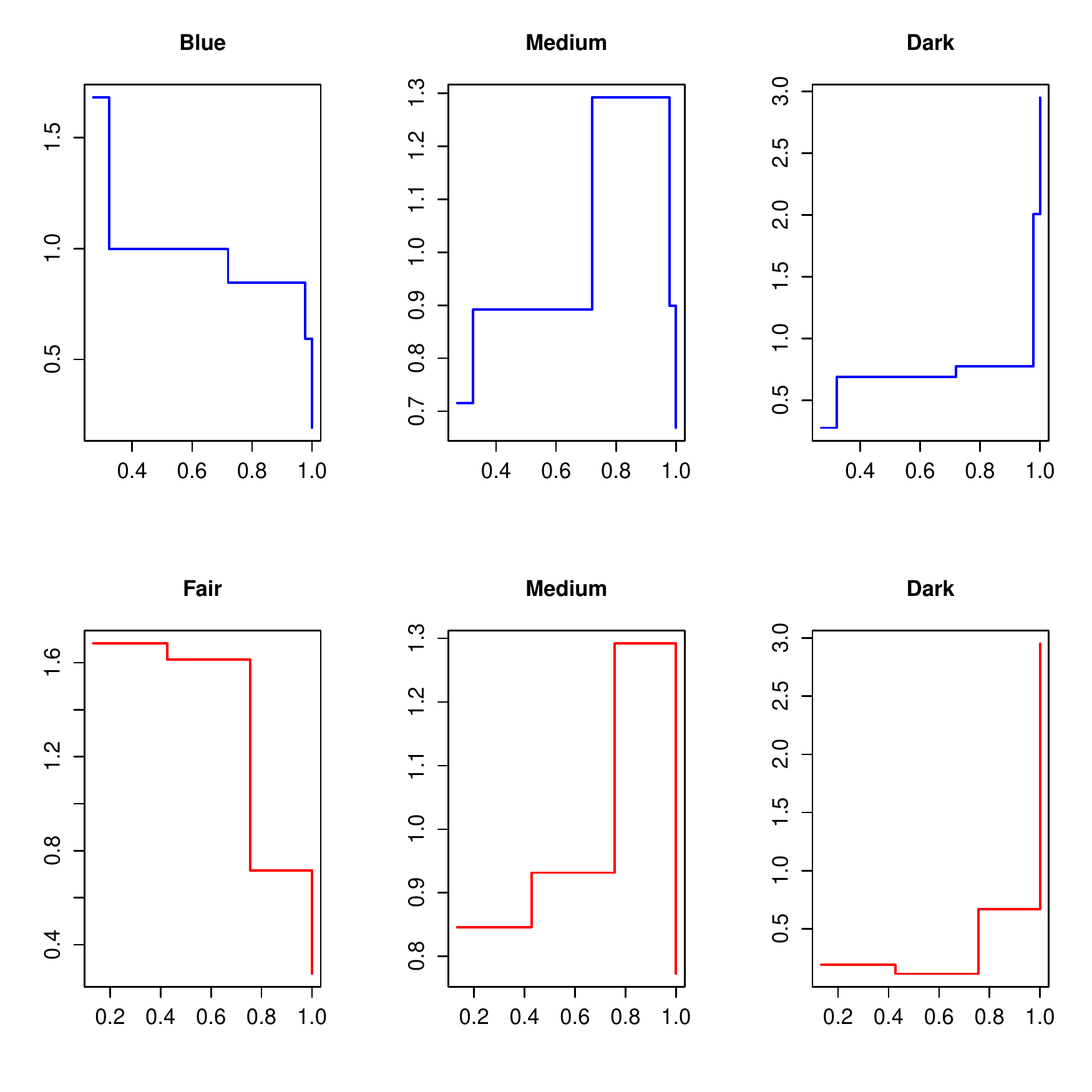}
\vspace{-.1em}
\caption{(color online) Correspondence analysis matches the shapes of conditional comparison density of row categories (blue) with column categories (red). See the similarity of the shapes of row and comumn categoriess and matches with the findings in Fig \ref{fig:fcorsp} } \label{fig:fccd}
\end{figure}

Orthogonal $L^2$ coefficients of the piece-wise constant conditional copula densities (slices of copula density) of our canonical copula model satisfies:
\bea
\LP[k;Y,Y|X=Q(u;X)]&=& \la_k \,\phi_k(u;X),\\
\LP[k;X,X|Y=Q(v;Y)]&=& \la_k \,\psi_k(v;X).
\eea
This implies that the LP correspondence scores capture the shape of  piece-wise constant conditional copula density functions. Thus, our correspondence analysis algorithm can be interpreted as comparing the `similarity of the shapes of conditional comparison density functions' (captured by orthogonal coefficients) of row categories with column categories, as shown in the Fig \ref{fig:fccd} for Fisher data.

\subsubsection{Low-Rank Smoothing Model, Smart Computational Algorithm}
All currently available methods perform singular value decomposition on \emph{raw} distance matrix (also known as contingency ratios) $\cop(F(x),F(y)) = p(x,y;X,Y)/p(x;X)p(y;Y)$ or $\log[\cop(F(x),F(y))]$, which is of the order $I \times J$;  this can have prohibitively high computational complexity and storage requirements for very large $I$ and $J$. Moreover, this is numerically unstable for sparse tables.
\vskip.25em

On the other hand, LP algorithm performs correspondence analysis by SVD on LP comoment matrix, which is often a \emph{far smaller order than the dimensional of the observed frequency matrix}. Thus, our algorithm provides efficient ways of compressing the structure of large tables. For the structured contingency table, the spectrum of ``raw'' contingency ratio matrix can be approximated by the spectrum of ``smooth'' LP-comoment matrix, which could provide enormous computational gain for large data.

\subsection{Interpretation: $2 \times 2$ Table }
We establish connection between the parameters of our LP copula models and traditional statistical measures for $2 \times 2$ table. Our approach provides an alternative elegant derivation presented in \cite{gilula1988} and \cite{goodman1991}.
\begin{thm} For $2 \times 2$ contingency table with Pearson correlation $\phi$ and odds-ratio $\de$,
\bea
\la_1&=&\big|\LP[1,1;X,Y]\big|\,\,=\,\,\big|\phi(X,Y)\big| \\
\ga_1&=&\big|\te_{1,1}\big|\,\,=\,\, \big| \log \de \big| \big( P_{1+}P_{+1}P_{2+}P_{+2} \big)^{1/2}
\eea
\end{thm}
To prove the result first verify that
\bea
T_1(0;X)=-\sqrt{P_{2+}/P_{1+}}, & T_1(1;X)=\sqrt{P_{1+}/P_{2+}} \nonumber \\
T_1(0;Y)=-\sqrt{P_{+2}/P_{+1}}, & T_1(1;Y)=\sqrt{P_{+1}/P_{+2}},
\eea
where $P_{i+}=\sum_j P_{ij}$ and $P_{+j}=\sum_i P_{ij}$, and $P_{ij}$ denotes the $2 \times 2$ probability table. We then have the following,
\bea
\LP[1,1;X,Y]&=&\Ex[T_1(X;X) T_1(Y;Y)]\,=\,\sum_{i=1}^2 \sum_{j=1}^2 P_{ij} T_1(i-1;X)\,T_1(j-1;Y) \nonumber \\
&=& \Big[P_{11}(P_{2+}P_{+2}) - P_{12}(P_{2+}P_{+1}) - P_{21}(P_{1+}P_{+2}) + P_{22}(P_{1+}P_{+1})   \Big]/\big( P_{1+}P_{+1}P_{2+}P_{+2} \big)^{1/2}\nonumber \\[.15em]
&=&\big( P_{11}P_{22}-P_{12}P_{21}\big)/\big( P_{1+}P_{+1}P_{2+}P_{+2} \big)^{1/2}\,=\,\phi(X,Y).
\eea
To prove (5.6) verify that
\bea
\te_{1,1}&=&\int \log\big[\cop(u,v)\big] S_1(u;X)S_1(v;Y)\dd u \dd v \nonumber \\
&=& \sum_{i=1}^2 \sum_{j=1}^2 \log \left[\dfrac{P_{ij}}{P_{i+}P_{+j}}\right] \, T_1(i-1;X)\,T_1(j-1;Y) \,P_{i+}\,P_{+j} \nonumber \\
&=&\log \left[\dfrac{P_{11}P_{22}}{P_{12}P_{21}}\right]\,  \big( P_{1+}P_{+1}P_{2+}P_{+2} \big)^{1/2}.
\eea
\section{LPINFOR Dependence Measure}
Computationally fast copula-based nonparametric dependence measure LPINFOR and its properties are discussed in this section. We present two decompositions of LPINFOR, which allow finer understanding of the dependence between $X$ and $Y$ from two different perspectives.
Applications to nonlinear dependence modeling and sparse contingency table are considered. A comparative simulation study is also
provided.
\subsection{Definition, Properties and Interpretation}
Dependence of $X$ and $Y$ can be nonparametrically measured by
\beq \label{eq:lpinfor}
\LPINFOR(X,Y)~=~\sum_{j,k > 0} \big|\LP[j,k;X,Y]\big|^2~=~\iint_{[0,1]^2} \big[\cop(u,v;X,Y) -1\big]^2 \dd u \dd v.\eeq
The second equality follows from LP copula representation \eqref{thm:lpcop} by applying Parseval's identity. Summing over significant indices $j,k$ gives a smooth estimator of $\LPINFOR(X,Y)$. The power of this method is best described by its applications, which are presented in the next section. The name LPINFOR is motivated from the observation that it can be interpreted as an INFORmation-theoretic measure belonging to the family of Csiszar's f-divergence measures \citep{csiszar1975}. LPINFOR is simultaneously a global correlation index and a directional statistic whose components capture the ways in which the copula density deviates from uniformity. By construction, LPINFOR is invariant to monotonic transformations of the two variables.

An alternative representation of LPINFOR can be derived from the singular value decomposition representation (5.1) as $\LPINFOR(X,Y)=\sum_k \la^2_k$. Note that our measure works for discrete or continuous marginals. LPINFOR, under two important cases (bivariate normal and contingency table), are described in the next result.

\begin{thm}[Properties]\label{thm:lpp} Two important properties for LPINFOR measure
\vskip.15em
(i) For $(X,Y)$ bivariate normal $\LPINFOR$ is an increasing function of $|\rho|$:
\beq \LPINFOR(X,Y)\, =\, \dfrac{\rho^2}{1-\rho^2}. \eeq
\vskip.15em
(ii) For $(X,Y)$ forming a $I \times J$ contingency table we have
\beq \LPINFOR(X,Y)\,=\,\CHIDIV(X,Y)\,=\,\chi^2/n.\eeq
\end{thm}
Part (i): To compute the integral $\int \big[\cop(u,v) - 1\big]^2 du\,dv$ for Gaussian copula, a quantity that is invariant to the choice of orthonormal polynomials, we first switch orthogonal polynomials from LP- orthonormal score functions to Hermite polynomials and use the following well-known result to complete the proof.

\begin{lemma}[Eagleson, 1964\nocite{eagleson64} and Koziol, 1979\nocite{koziol1979}]\label{lem:gcop} Let $H_j(x)$ denotes the $j$th orthonormalized Hermite polynomial. For bivariate Gaussian copula with correlation parameter $\rho$ we have the following inner-product result:
\[ \iint_{[0,1]^2} \dfrac{\phi[ \Phi^{-1}(u;X),\Phi^{-1}(v;Y); \rho   ]}{\phi[\Phi^{-1}(u;X)]\, \phi[\Phi^{-1}(v;Y)]}\,H_j[\Phi^{-1}(u;X)]\,H_k[\Phi^{-1}(v;Y)] \dd u \dd v\,=\,\rho^{j}\, \ind[j=k]. \]
\end{lemma}
\vskip.25em
Part(ii): Note that for $(X,Y)$ discrete, $\cop(F(x),F(y);X,Y)=p(x,y;X,Y)/p(x;Y)p(y;Y)$. Verify the following representation to complete the proof
\beq \CHIDIV(X,Y)\,=\,\chi^2/n\,=\,\sum_{x,y}\, p(x;Y) \, p(y;Y) \big[\cop\big(F(x),F(y);X,Y\big)  - 1\big]^2.  \eeq

\vskip.5em

{\bf Example.} For the Ripley data we previously computed the LP-comoment matrix (2.8). The smooth LPINFOR (taking sum of squares of the significant elements denoted by `$^*$') dependence number is $1.851$. There is strong evidence of nonlinearity due to the presence of higher-order LP-comoments. Coefficient of linearity can be defined as $|\LP[1,1]|^2/\LPINFOR = .45$. Other applications of LPINFOR as a nonlinear correlation measure will be investigated in Section 6.3.2.
\vskip.25em

For Fisher data, the chisquare divergence or total inertia: $\chi^2/n$=.230 with df $(4-1) \times (5-1)=12 $, and the smooth-LPINFOR dependence number is $.220$, with \emph{much fewer degrees of freedom} $3$ (2.9). LPINFOR numerically computes the `raw' chisquare  divergence number but the computation of degrees of freedom is different. In our case, degrees of freedom dramatically drops from the usual $12$ to $3$, which boosts the power of our method. Applicability for large sparse tables is described in Section 6.3.1.

\subsection{Conditional LPINFOR, Interpretation}
Define conditional LPINFOR
\beq \LPINFOR(Y|X=x)\,=\,\int_0^1 \big[ d\left( v;\,Y,Y|X=x\right)\,-\, 1\big]^2 \dd v \,=\,\sum_{k>0} \big|\LP[k;Y,Y|X=x]\big|^2.\eeq
The components $\LP[k;Y,Y|X=x]=\Ex[T_k(Y;Y);Y|X=x]$ computed from LP-comoments by \[\LP[k;Y,Y|X=x]~=~\sum_{j>0} \LP[j,k;X,Y]\,T_j(x;X).\]
The LPINFOR can be decomposed as
\beq
\LPINFOR(X,Y)~=~\int  \LPINFOR(Y|X=x)\dd F(x;X)
\eeq

The conditional LPINFOR quantifies the effect of different levels of $X$ on $Y$. For $(X,Y)$ two-way continency table the components of conditional LPINFOR $\LP[k;Y,Y|X=x]$ provide the LP profile coordinates for correspondence analysis (Section 5.2).
For $(X,Y)$ continuous these components give insight how the shape of conditional density $f(y;Y|X=x)$ changes with the levels of the covariate $X$, illustrated in Section 7.4.
\subsection{Applications}
\subsubsection{Sparse Contingency Table} \label{sec:lst}
Consider the example shown in Table \ref{tab:iqt}.  Our goal is to understand the association between Age and Intelligence Quotient (IQ). It is well known that the accuracy of chi-squared approximation depends on both sample size $n$ and number of cells $p$.
Classical association mining tools for contingency tables like Pearson chi-squared statistic and likelihood ratio statistic break down for sparse tables (where many cells have observed zero counts, see \cite{haberman88})  - a prototypical example of ``large $p$ small $n$'' problem.

The classical chi-square statistic yields test statistic = $60$, df = $56$, p-value = $0.3329$, failing to capture the pattern from the sparse table. Eq \ref{eq:lpiq} computes the LP-comoment matrix. The significant element is $\LP[2,1]=-0.617$, with the associated p-value $0.0167$. LPINFOR thus adapts successfully to sparse tables.

\begin{table}[tbh!]
\setlength{\tabcolsep}{5.5pt}
\caption{Wechsler Adult Intelligence Scale (WAIS) by males of various ages. Is there any pattern ? Data described in  \cite{mack1981} and \cite{rayner96}.}

\begin{center}
\begin{tabular}{cccccccccccccccc}
\toprule
 &\multicolumn{15}{c}{Intelligence Scale}\\
\cmidrule(r){2-16}
Age group &1 & 2 &3 &4&5&6&7&8&9&10&11&12&13&14&15\\
\hline
$16 - 19$ & 0 &1 &0 &0 &0 &0 &0 &1 &0 &1 &0 &0 &0 &0 &0\\
$20 - 34$ & 0 &0 &0 &0 &0 &1 &0 &0 &0 &0 &1 &0 &0 &1 &0\\
$35 - 54$ &  0 &0 &0 &0 &0 &0 &0 &0 &1 &0 &0 &0 &1 &0 &1\\
$55 - 69$ & 0 &0 &1 &0 &0 &0 &1 &0 &0 &0 &0 &1 &0 &0 &0\\
$69 +$    & 1 &0 &0 &1 &1 &0 &0 &0 &0 &0 &0 &0 &0 &0 &0\\
\hline
\end{tabular}
\end{center}
\label{tab:iqt}
\end{table}

\beq \label{eq:lpiq}
\widehat{\LP}\big[X={\rm Age} ,Y={\rm IQ}  \big]~=~\begin{bmatrix}
                               -0.316  &0.173  &0.168 &-0.114 \\[.38em]
                              -0.618^{*} &-0.031 &-0.101  &0.068\\[.38em]
                               0.087  &0.136  &0.077  &0.037\\[.38em]
                               0.165  &0.215  &0.042  &0.289\\[.38em]
                             \end{bmatrix}\eeq

The higher-order significant $\LP[2,1;X,Y]=\Ex[T_2(X;X)T_1(Y;Y)]$ comoment not only indicates strong evidence of non-linearity but also gives a great deal of information regarding the functional relationship between age and IQ. Our analysis suggests that IQ is associated with the $T_2(X;X)$ that has a quadratic (more precisely, inverted U-shaped) age effect, which agrees with all previous findings that IQ increases as a function of age up to a certain point, and then it declines with increasing age. This example also demonstrates how the LP-comoment based approach  nonparametrically finds the optimal non-linear transformations for two-way contingency tables, extending \cite{breiman85}
  \[ \rho\big( \Psi^{*}(Y),\, \Upsilon^{*}(X)\big)~=~\max_{\Psi,\phi} \rho \big( \Psi(Y),\, \Upsilon(X)\big). \]

\subsubsection{Nonlinear Dependence Measure, Power Comparison} \label{sec:simu}
\begin{figure}[tbh!]
 \centering
 \includegraphics[height=.45\textheight,width=.85\textwidth,trim=1.5cm .5cm 1.5cm .5cm]{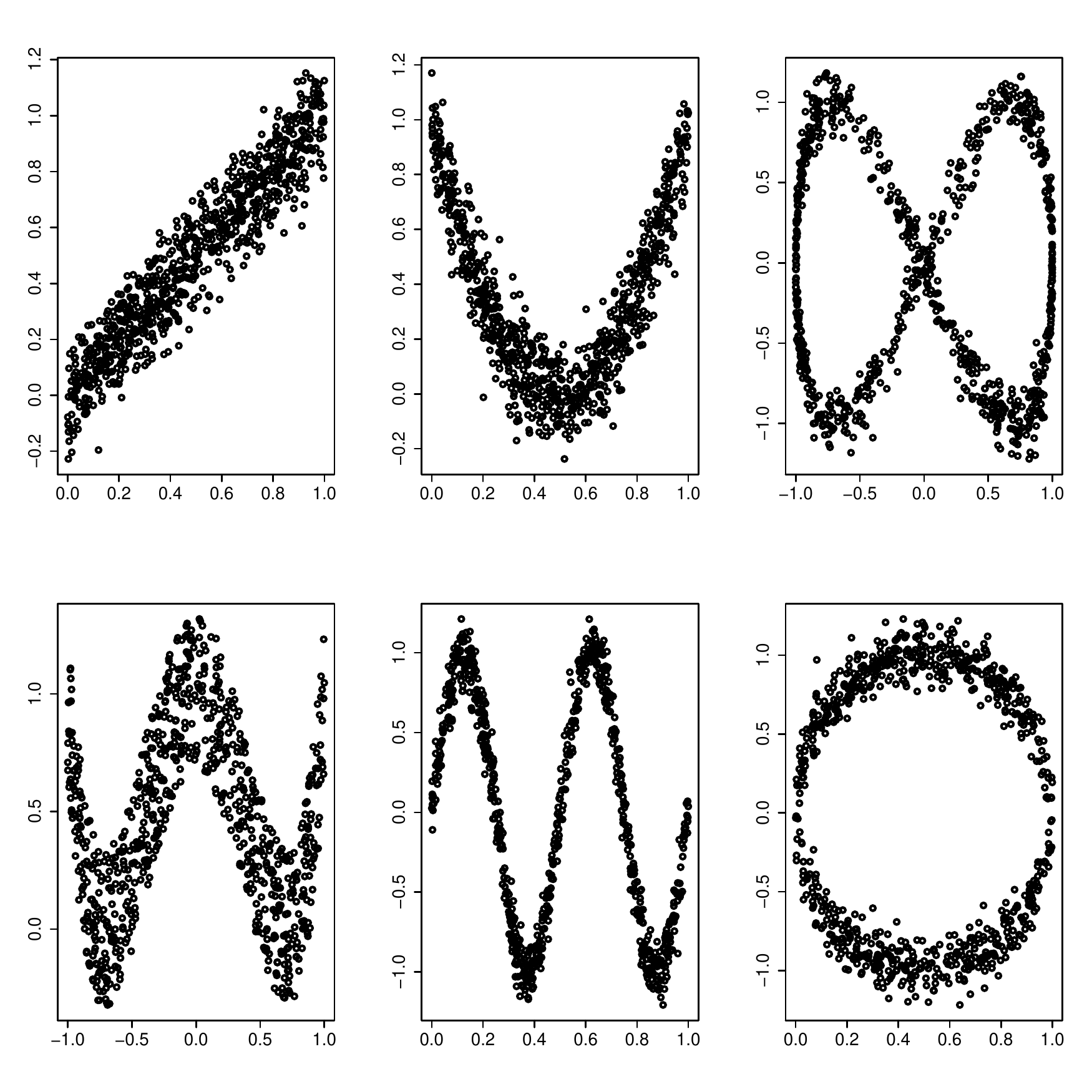}
\vspace{-.1em}
\caption{Various functional relationships considered in our simulation study in Section \ref{sec:simu}. }
\label{fig:func}
\end{figure}

We compare the power of our LPINFOR statistic with distance correlation \citep{szekely09}, maximal information coefficient \citep{MI11}, Spearman and Pearson correlation for testing independence under six different functional relationships (for brevity we restrict our attention to the relationships shown in Fig \ref{fig:func}) and four different types of noise settings: (i) $E_1:$ Gaussian noise $\cN(0,\si)$ where $\si$ varying from $0$ to $3$; (ii) $E_2:$ Contaminated Gaussian $(1-\eta)\cN(0,1)+\eta \cN(1,3)$, where the contamination parameter $\eta$ varies from $0$ to $.4$; (iii) $E_3:$ Bad leverage points are introduced to the data by $(1-\eta)\cN(0,1)+\eta \cN(\mu,1)$, where $\mu=\{\pm20,\pm40\}$ w.p $1/4$ and $\eta$ varies from $0$ to $.4$; (iv) $E_4:$ heavy-tailed noise is introduced via Cauchy distribution with scale parameter that varies from $0$ to $2$.

\begin{table}[thb!]
\setlength{\tabcolsep}{4.5pt}
\caption{Performance summary table on the basis of empirical power among six competing methods, compared over six bivariate relationships and four different noise settings.}
\vskip1.5em
\centering
\begin{tabular}{c c cccccc}
\toprule
\toprule
 & &\multicolumn{6}{c}{Functional Relation}\\
\cmidrule(r){3-8}
Noise & Performance &{\bf Linear} &{\bf Quadratic}  &{\bf Lissajous} &{\bf W-shaped}&{\bf Sine} &{\bf Circle}\\
\midrule
$E_1$ & Winner &Pearson &LPINFOR &LPINFOR & Dcor  &Dcor  & LPINFOR\\[1ex]
  & Runner-up  &Spearman &Dcor   &MIC     &LPINFOR &MIC & MIC \\[1ex]
\hline\\ [-1.5ex]
$E_2$ & Winner &Spearman &LPINFOR &LPINFOR & LPINFOR  &MIC  & LPINFOR\\[1ex]
  & Runner-up  &Dcor     &MIC      &MIC    & Dcor     &Dcor & MIC \\[1ex]
\hline\\ [-1.5ex]
$E_3$ & Winner &Spearman &LPINFOR &LPINFOR & LPINFOR  &MIC  & LPINFOR\\[1ex]
  & Runner-up  &LPINFOR     &MIC      &MIC    & MIC    &LPINFOR & MIC \\[1ex]
  \hline\\ [-1.5ex]
$E_4$ & Winner &Spearman &LPINFOR &LPINFOR & LPINFOR  &MIC  & LPINFOR\\[1ex]
  & Runner-up  &LPINFOR     &MIC      &Dcor    & MIC    &LPINFOR & MIC \\[1ex]
   \hline\\ [-1.5ex]
\end{tabular}
\label{tab:power1}
\vskip.5em
\end{table}

The simulation is carried out based on sample size $n=300$. We used $250$ null simulated data sets to estimate the 95\% rejection cutoffs for all methods at the significance level $0.05$ and used $200$ simulated data sets from alternative to approximate the power, similar to \cite{Tibs11}. We have fixed $m=4$ to compute LPINFOR for all of our simulations. Fig 10-12 show the power curves for the bivariate relations under four different noise settings. We summarize the findings in Table \ref{tab:power1}, which provides strong evidence that LPINFOR performs remarkably well across s wide-range of settings.

\vskip.5em
{\bf Computational Considerations.} We investigate the issue of scalability to large data sets (big $n$ problems) for the three competing nonlinear correlation detectors: LPINFOR, distance correlation (Dcor), and maximal information coefficient (MIC).
Generate $(X,Y) \sim \mbox{Uniform}[0,1]^2$ with sample size $n=20,50,100,500,1000,5000$ $,10,000$. Table \ref{tab:time} reports the average time in sec (over $50$ runs) required to compute the corresponding statistic. It is important to note the software platforms used for the three methods: DCor is written in JAVA, MIC is written in C, and we have used (unoptimized) R for LPINFOR. It is evident that for large sample sizes both Dcor and MIC become almost infeasible - not to mention computationally extremely expensive. On the other hand, LPINFOR emerges clearly as the computationally most efficient and scalable approach for large data sets.

\begin{table}[t!]
\setlength{\tabcolsep}{5pt}
\caption{Computational Complexity: uniformly distributed, independent samples of size $n$, averaged over $50$ runs based on Intel(R) Core(TM) i7-3540M CPU @ 3.00GHz 2 Core(s) processor. Timings are reported in seconds.}
\vskip1.5em
\centering
\begin{tabular}{c c cccccc}
\toprule
\toprule
Methods&\multicolumn{6}{c}{Size of the data sets}\\
\cmidrule(r){2-7}
&$n=100$ & $n=500$ & $n=1000$ & $n=2500$ & $n=5000$ & $n=10,000$\\
\midrule
LPINFOR & 0.002 (.004) & 0.001 (.004) & 0.002 (.005) & 0.005 (.007)  & 0.011 (.007)  &0.018 (.007)  \\[1ex]

Dcor  & 0.003 (.007) &0.094 (.013) &0.371 (.013) &1.773 (.450)  &7.756 (.810) &44.960 (12.64) \\[1ex]

MIC & 0.312 (.007) & 0.628 (.008) & 1.357 (.035) & 5.584 (.110) & 19.052 (.645) & 65.170 (2.54)   \\[1ex]
\hline
\end{tabular}
\label{tab:time}
\vskip.5em
\end{table}

\begin{figure}[ht]
 \centering
 \vspace{-.15em}
 \subfigure[Power curve for linear pattern]{
  \includegraphics[height=.44\textheight,width=.85\textwidth,trim=1cm 0cm 1cm 1cm]{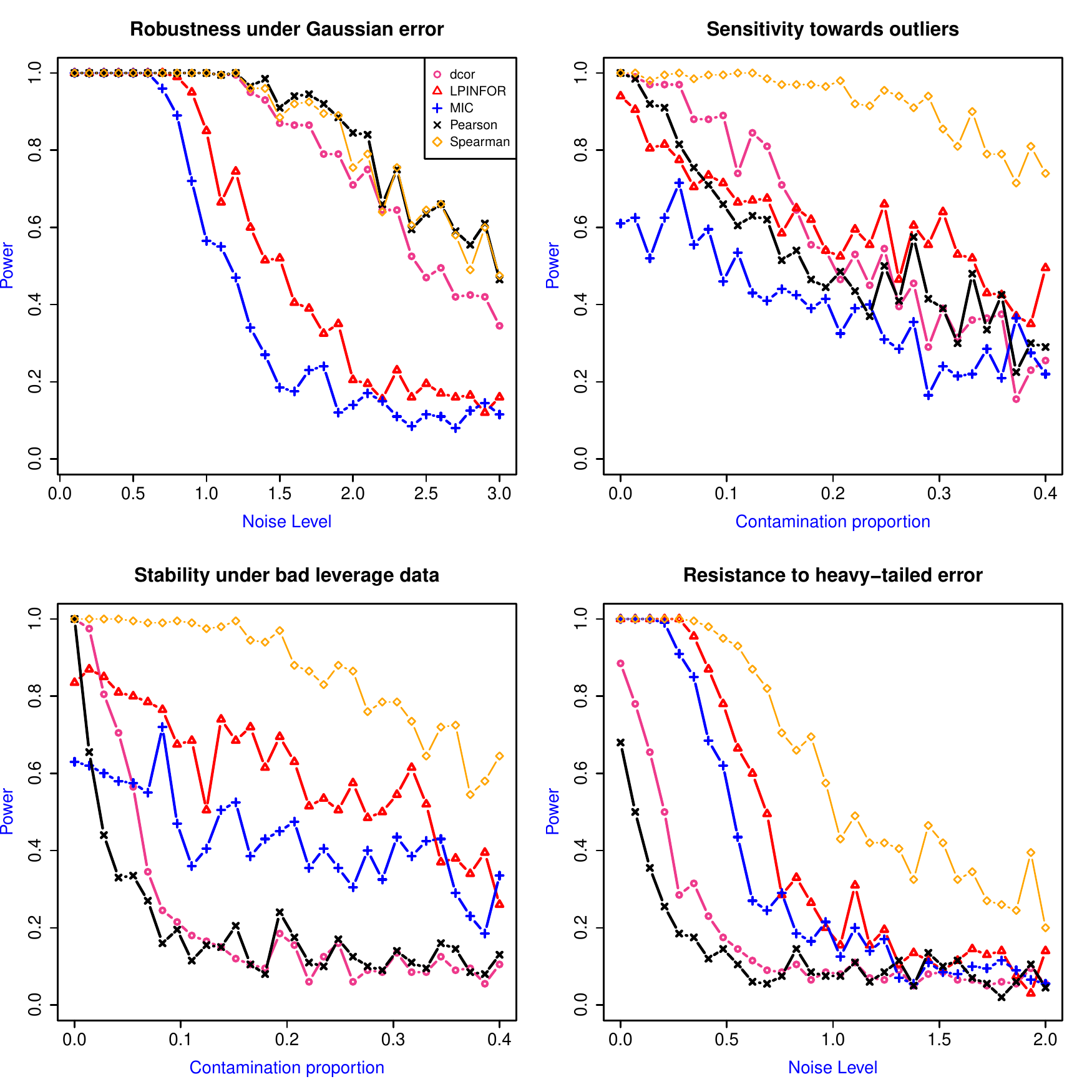}
   \label{fig:L}
   }
 \subfigure[Power curve for quadratic pattern]{
  \includegraphics[height=.44\textheight,width=.85\textwidth,trim=1cm .5cm 1cm 0cm]{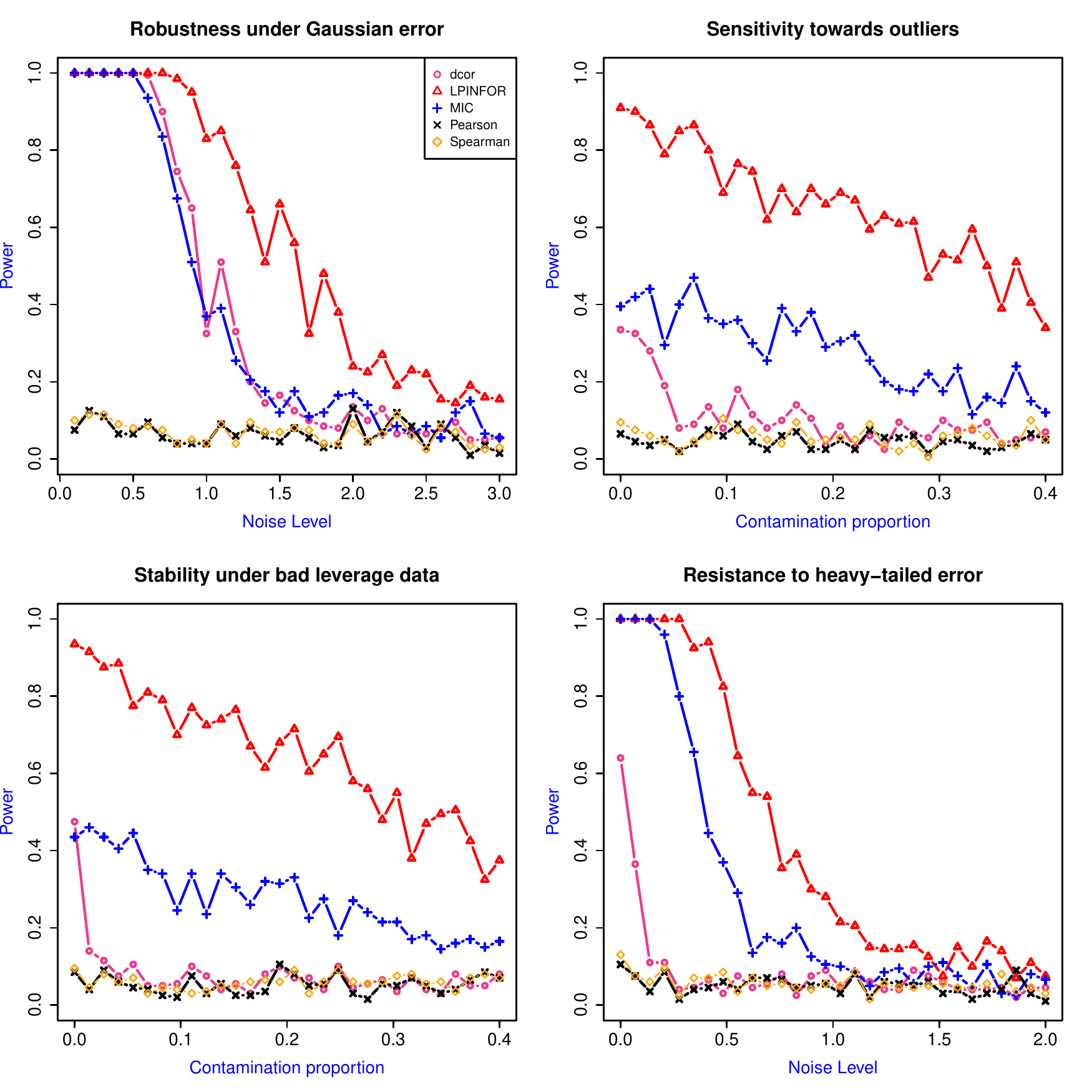}
   \label{fig:Q}
   }
 \caption{Estimated power curve for linear and quadratic pattern under four types of noise settings; set up described in Section 6.3.2}
\end{figure}

\begin{figure}[ht]
 \centering
 \vspace{-.15em}
  \subfigure[Power curve for Lissajous curve]{
  \includegraphics[height=.44\textheight,width=.85\textwidth,trim=1cm .5cm 1cm 0cm]{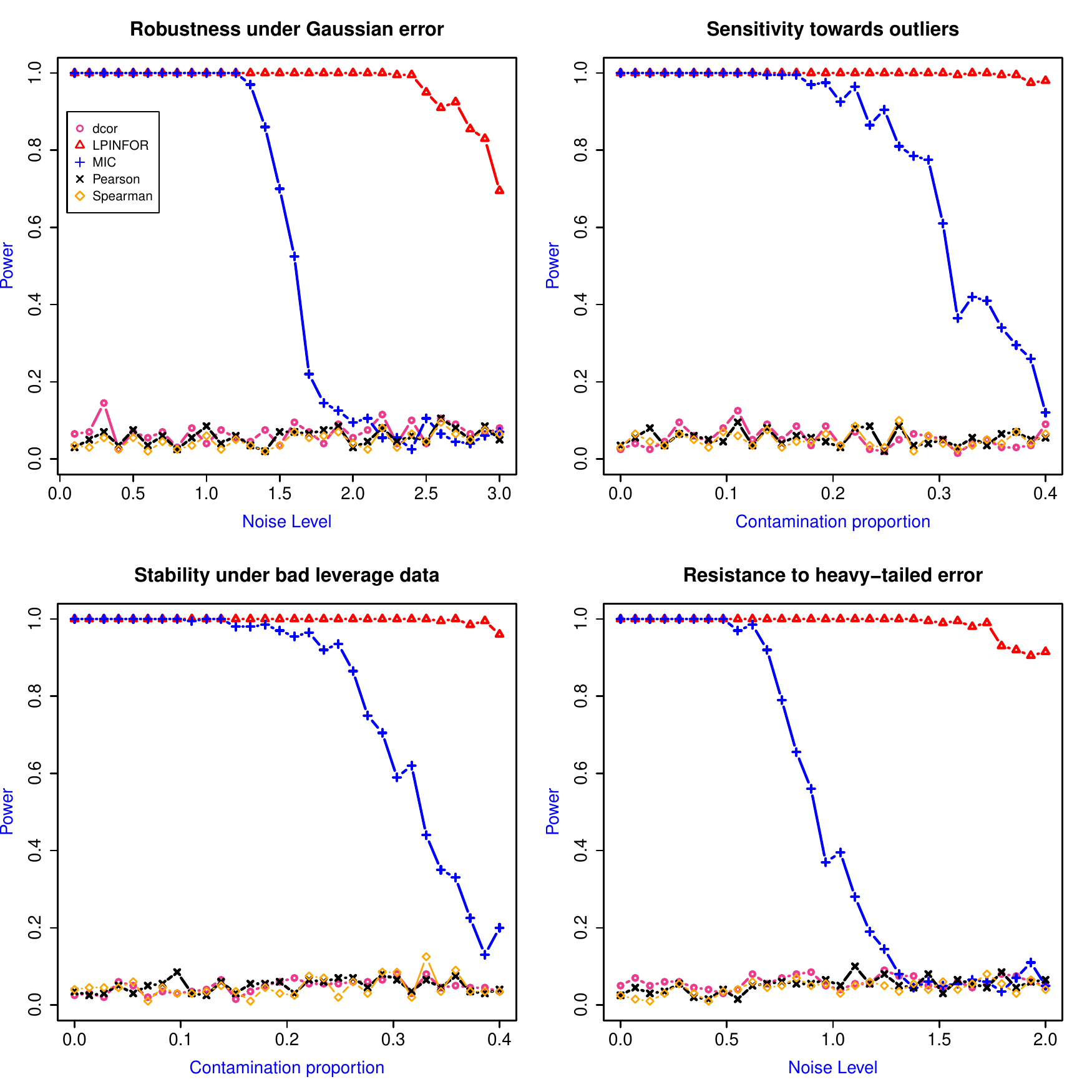}
   \label{fig:C}
   }
 \subfigure[Power curve for W-pattern]{
  \includegraphics[height=.44\textheight,width=.85\textwidth,trim=1cm 0cm 1cm 0cm]{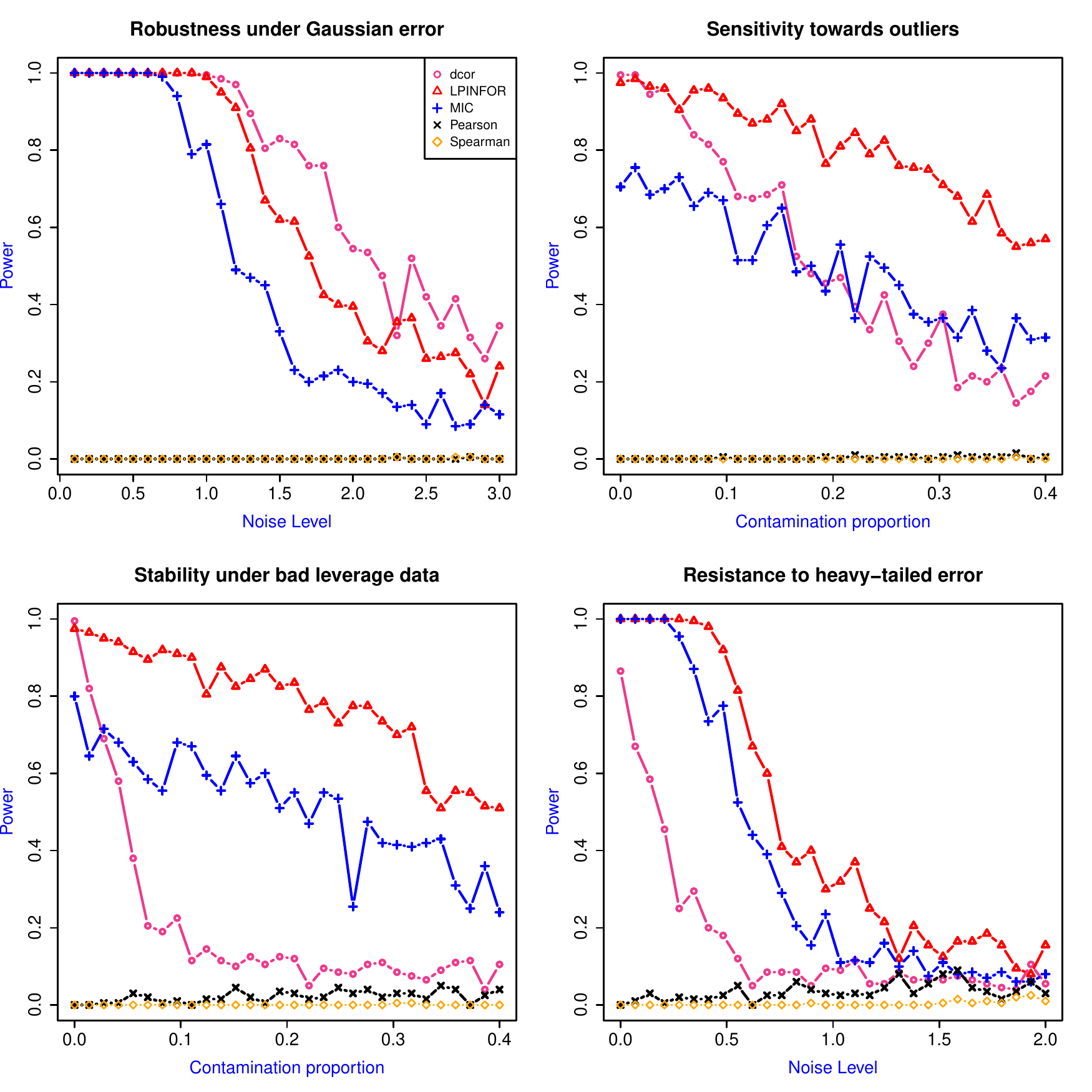}
   \label{fig:T}
   }
 \caption[]{Estimated power curve for Lissajous curve and W-pattern under four types of noise settings; set up described in  Section 6.3.2  }
\end{figure}

\begin{figure}[ht]
 \centering
 \vspace{-.15em}
 \subfigure[Power curve for sine pattern]{
  \includegraphics[height=.44\textheight,width=.85\textwidth,trim=1cm 0cm 1cm 1cm]{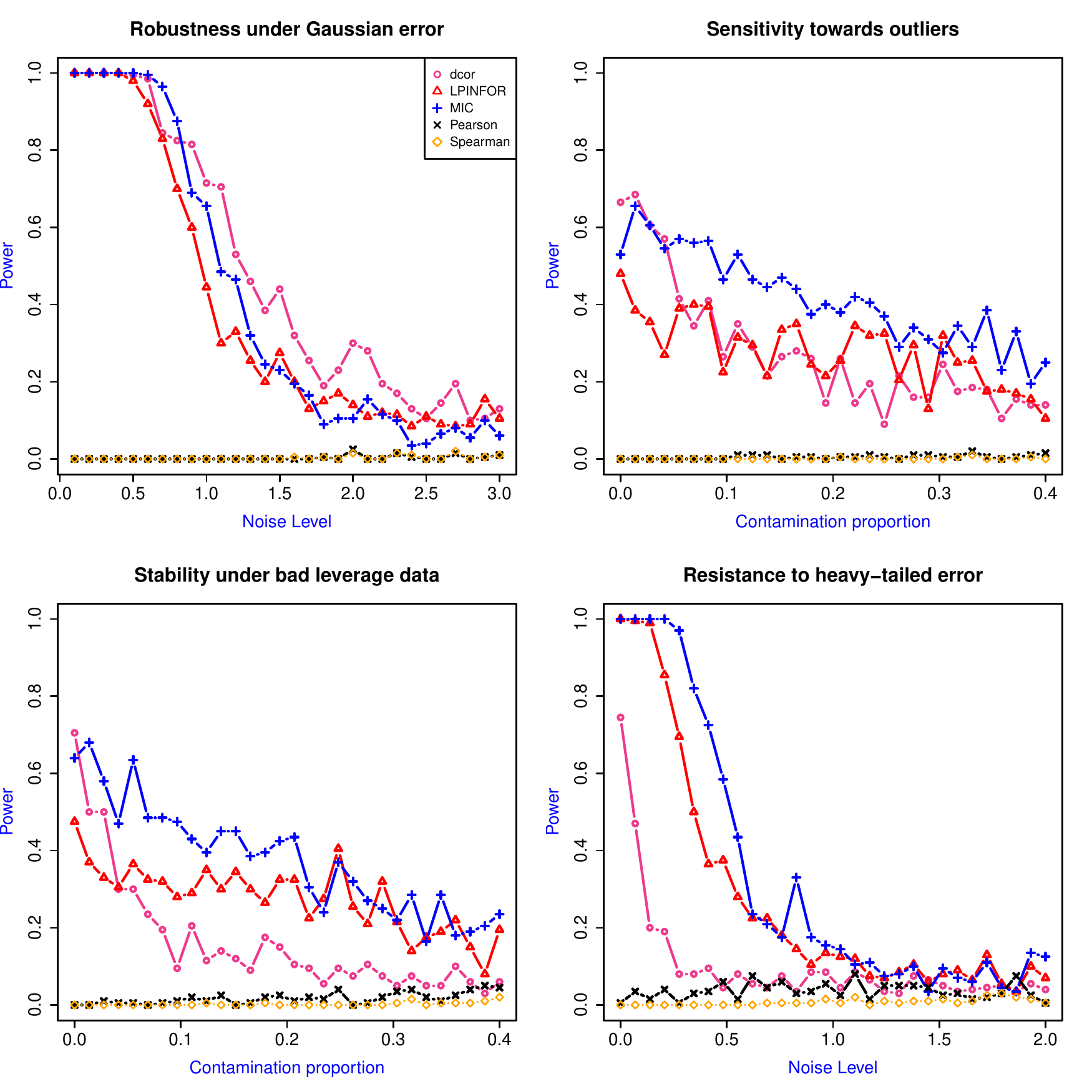}
   \label{fig:sine}
   }
 \subfigure[Power curve for circle pattern]{
  \includegraphics[height=.44\textheight,width=.85\textwidth,trim=1cm .5cm 1cm 0cm]{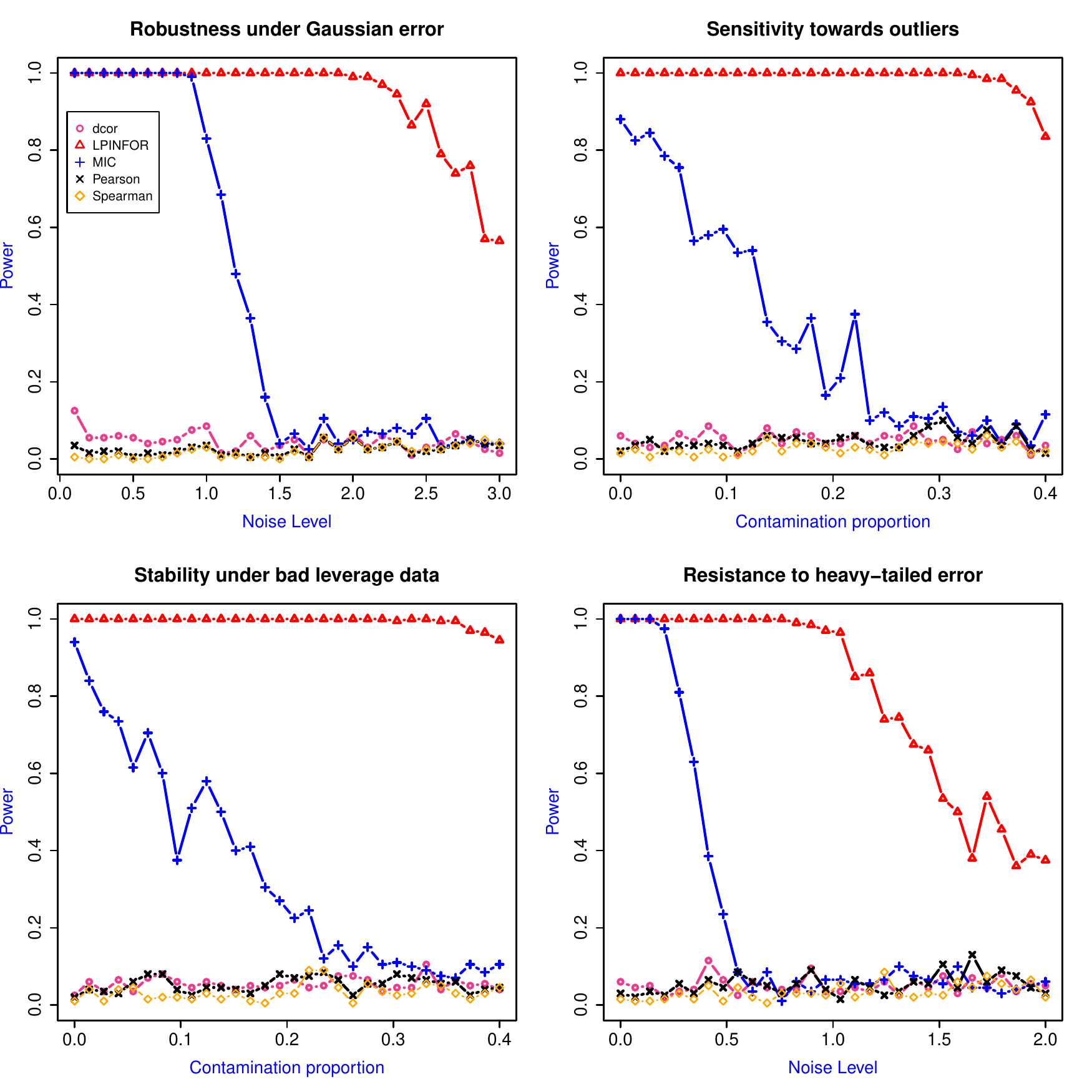}
   \label{fig:cir}
   }
 \label{fig:subfigureExample}
 \caption{Estimated power curve for sine curve and circle pattern under four types of noise settings; set up described in  Section 6.3.2}
\end{figure}

\section{LPSmoothing (X,Y)}

We propose a nonparametric copula-based conditional mean and conditional quantile modeling algorithm. We introduce LP zero-order comoments (defined for discrete and continuous marginals) as the regression coefficients of our nonparametric model (generalizing Gini correlation \cite{gini87} and \cite{Serf07}). An illustration using Ripley data is also given, demonstrating the versatility of our approach.

\subsection{Copula-Based Nonparametric Regression}
\begin{thm} \label{thm:nreg}
Conditional expectation $\Ex[Y|X=x]$ can be approximated by linear combination of orthonormal LP-score functions $\{T_j(X;X)\}_{1\le j \le m}$ with coefficients defined as $\Ex\big[\Ex[Y|X]T_j(X;X)\big] = \Ex[Y T_j(X;X)] =\LP(j,0;X,Y)$, because of the representation:
\beq \label{eq:cmean}
\Ex\big[ Y \mid X=x   \big] - \Ex[Y] \,=\, \sum_{j>0} T_j(x;X) \LP(j,0;X,Y).
\eeq
\end{thm}
Direct proof uses $\Ex\big[(Y-\Ex[Y|X]) T_j(X;X) \big]=0$ for all $j$, and Hilbert space theory.
An alternative derivation of this representation is given below, which elucidates the connection between copula-based nonlinear regression. Using LP-representation of copula density (4.3), we can express conditional expectation as
\bea
\Ex[Y|X=x]&=&\int y \cop\big( F(x), F(y) \big) \dd F(y) \,=\,\int y \Big[ 1+ \sum_{j,k>0}  \LP[j,k;X,Y]\, T_j(x;X)\, T_k(y;Y) \Big]  \dd F(y) \nonumber \\
&=& \Ex[Y]~+~\sum_{j>0} T_j(x;X)\, \left\{ \sum_{k>0} \LP[j,k;X,Y] \,\LP[k;Y]  \right\}.
\eea
Compete the proof by recognizing
\[\LP[j,0;X,Y]\,=\,\Ex\Big[\, \Ex\big[  T_j(X;X) | Y\big] Y \,\Big]\,=\,\sum_{k>0} \LP[j,k;X,Y] \,\LP[k;Y]. \]

Our approach is an alternative to parametric copula-based regression model \citep{noh2013,nikolo2010}. \cite{dette2013} aptly demonstrated the danger of the ``wrong'' specification of a parametric family even in one or two dimensional predictor. They noticed commonly used parametric copula families are not flexible enough to produce non-monotone regression function, which makes it unsuitable for any practical use.
Our LP nonparametric regression algorithm yields easy to compute closed form expression \eqref{eq:cmean} for conditional mean function, which bypasses the problem of misspecification bias and can adapt to any nonlinear shape. Moreover, our approach can easily be integrated with data-driven model section criterions like AIC or Mallow's $C_p$ to produce a parsimonious model.
\subsection{Generalized Gini Correlations as Zero Order LP-Comoments}
In this section We establish connection between the Zero-order LP-comoments (regression coefficients) and Gini correlation \citep{gini87} .

\emph{Generalized LP-Gini Correlations.} Define the Generalized Gini correlations by
\beas
\RGINI(j;Y|X) &=& \LP(j,0;X,Y)/\LP(j,0;Y,Y) = \Ex[T_j(X;X) Y]/\Ex[T_j(Y;Y)Y] \\
\RGINI(k;X|Y) &=& \LP(0,k;X,Y)/\LP(0,k;X,X) = \Ex[T_k(Y;Y) X]/\Ex[T_k(X;X)X]
\eeas
LP-comoment representation of traditional Gini correlation $\RGINI(1;Y|X)$ and $\RGINI(1;X|Y)$.
We extend the classical Gini correlation concept in at least two major directions: (i) our definition, unlike previous versions, is well-defined for both discrete and continuous variables. As noted in Chapter 2 of \cite{yitzhaki13book} to extend the concept of Gini for discrete random variables requires not-trivial adjustments, which are already incorporated in our LP-Gini algorithm. (ii) We provide non-linear generalizations (defined for $j \ge 1$) of classical Gini correlation using higher-order LP-score functions.

Next we investigate some useful properties of LP-Gini correlation for bivariate normal that can be used for quick diagnostic purpose.

\begin{thm} \label{thm:gini}
For $(X,Y;\rho)$ Bivariate Normal we have the following relationship
\[\LP[j,0;X,Y] \,=\, \rho\, \LP[j,0;Y,Y],~~\mbox{and}~~\LP[j,0;Y,X]\,=\, \rho\, \LP[j,0;X,X] .\]
\end{thm}
To prove the Theorem \ref{thm:gini} note that, for bivariate normal, $\Ex(Y|X) = \mu_y + \rho\, \si(y)\, X$ (w.l.g we can assume $X$ has mean $0$ and variance $1$ as it only appear as a function of ranks in Gini expression $\RGINI(j;Y|X)$).
For $Z\sim \cN(0,1)$ we then have
\beq \Ex[Y T_j(X)]= \Ex\big[\Ex[Y|X] T_j(X)\big] = \rho\, \si(y) \LP[0,j;Z,Z].\eeq
Finally, observe $\LP[j,0;Y,Y]= \si(y) \LP[0,j;Z,Z]$ to complete the proof.

\begin{coro} \label{cor:gini}
For $(X,Y;\rho)$ Bivariate Normal we have the following important property
\[\RGINI[j; Y|X]\,=\,\RGINI[j; X|Y]\,=\,\rho,~\mbox{for all}\,\, j\, \mbox{odd}.\]
\end{coro}

As a special case of Corollary \ref{cor:gini} we derive that the co-Ginis $\RGINI(1;Y|X)$ and $\RGINI(1;X|Y)$ are equal to $\rho$ for bivariate normal \citep{yitzhaki13book}.

\begin{figure}[ht]
 \centering
 \vspace{-.15em}
  \subfigure[Conditional Mean]{
  \includegraphics[height=.4\textheight,width=.45\textwidth,trim=1cm .5cm .5cm .5cm]{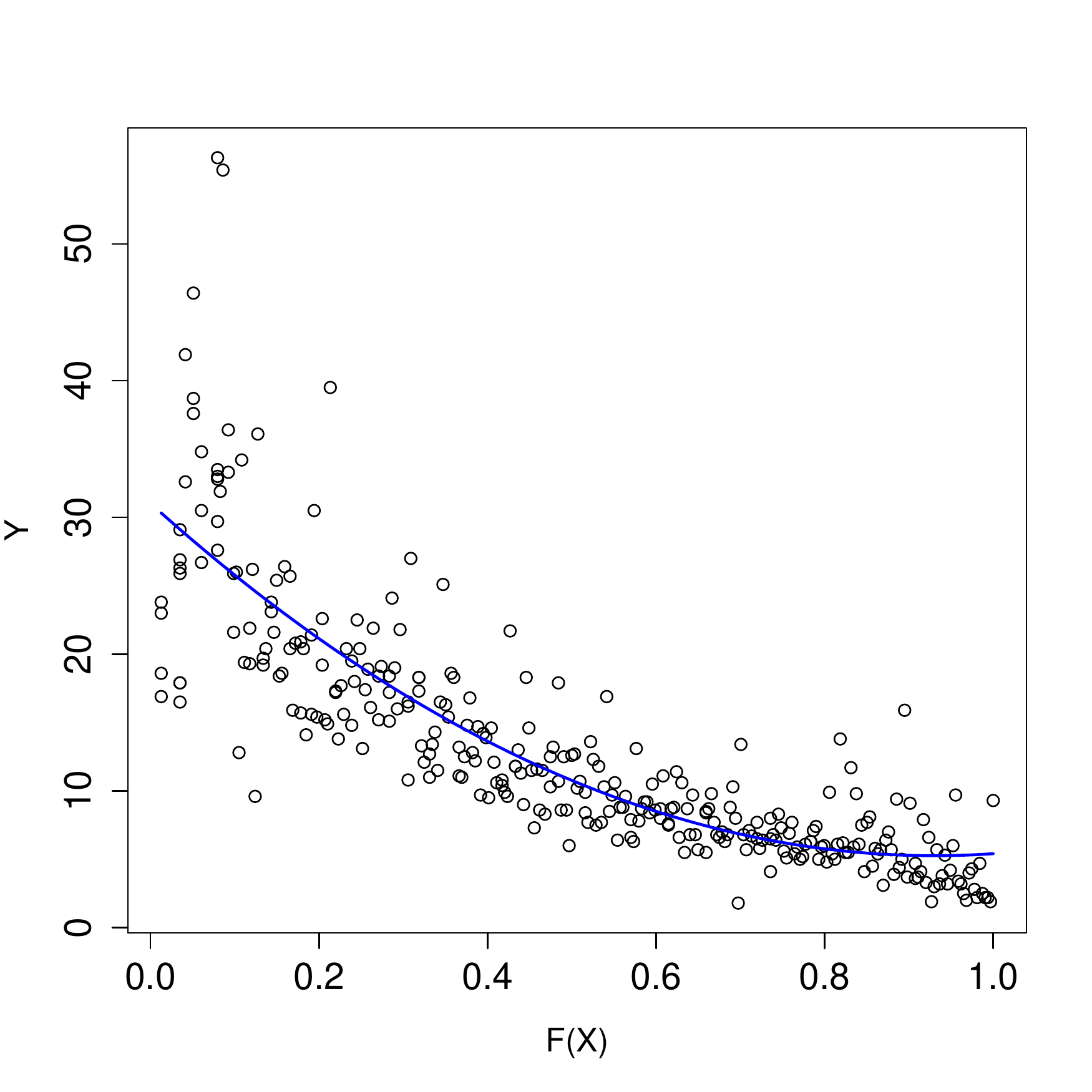}
   }
 \subfigure[Conditional Mean]{
  \includegraphics[height=.4\textheight,width=.45\textwidth,trim=.5cm .5cm .5cm .5cm]{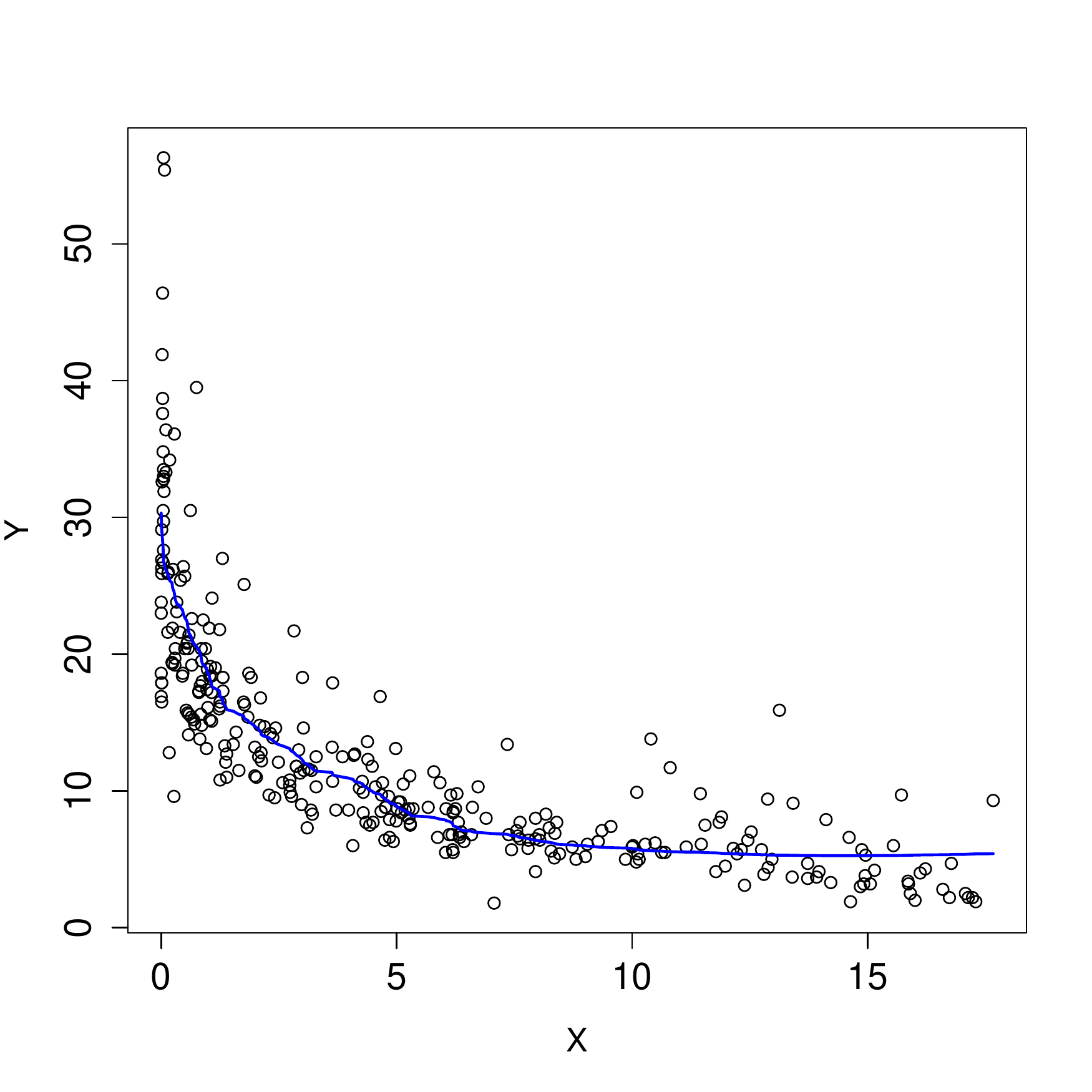}
   }\\

   \subfigure[Conditional Quantile]{
  \includegraphics[height=.4\textheight,width=.45\textwidth,trim=.5cm .5cm .5cm .5cm]{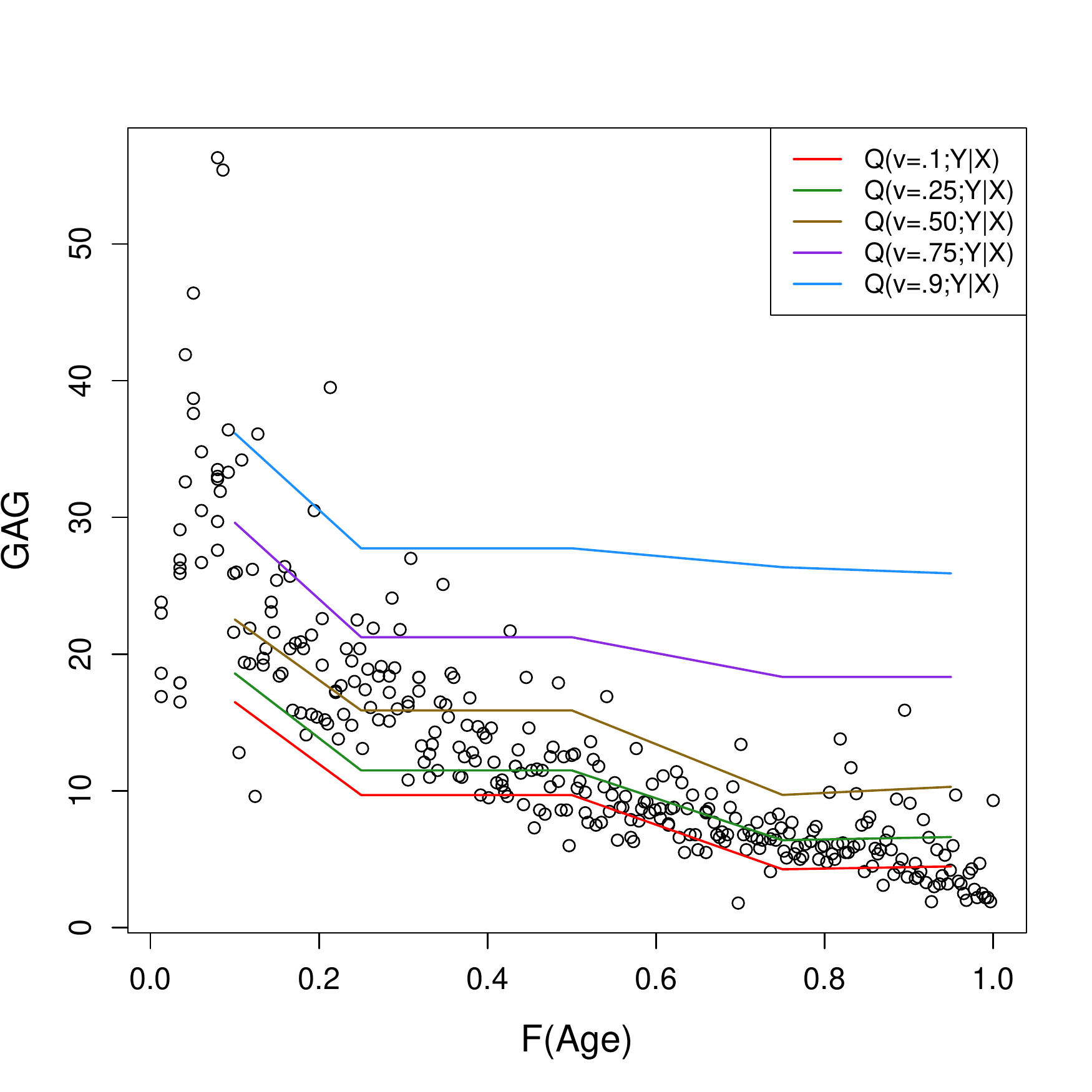}
   }
 \subfigure[Conditional Quantile]{
  \includegraphics[height=.4\textheight,width=.45\textwidth,trim=.5cm .5cm .5cm .5cm]{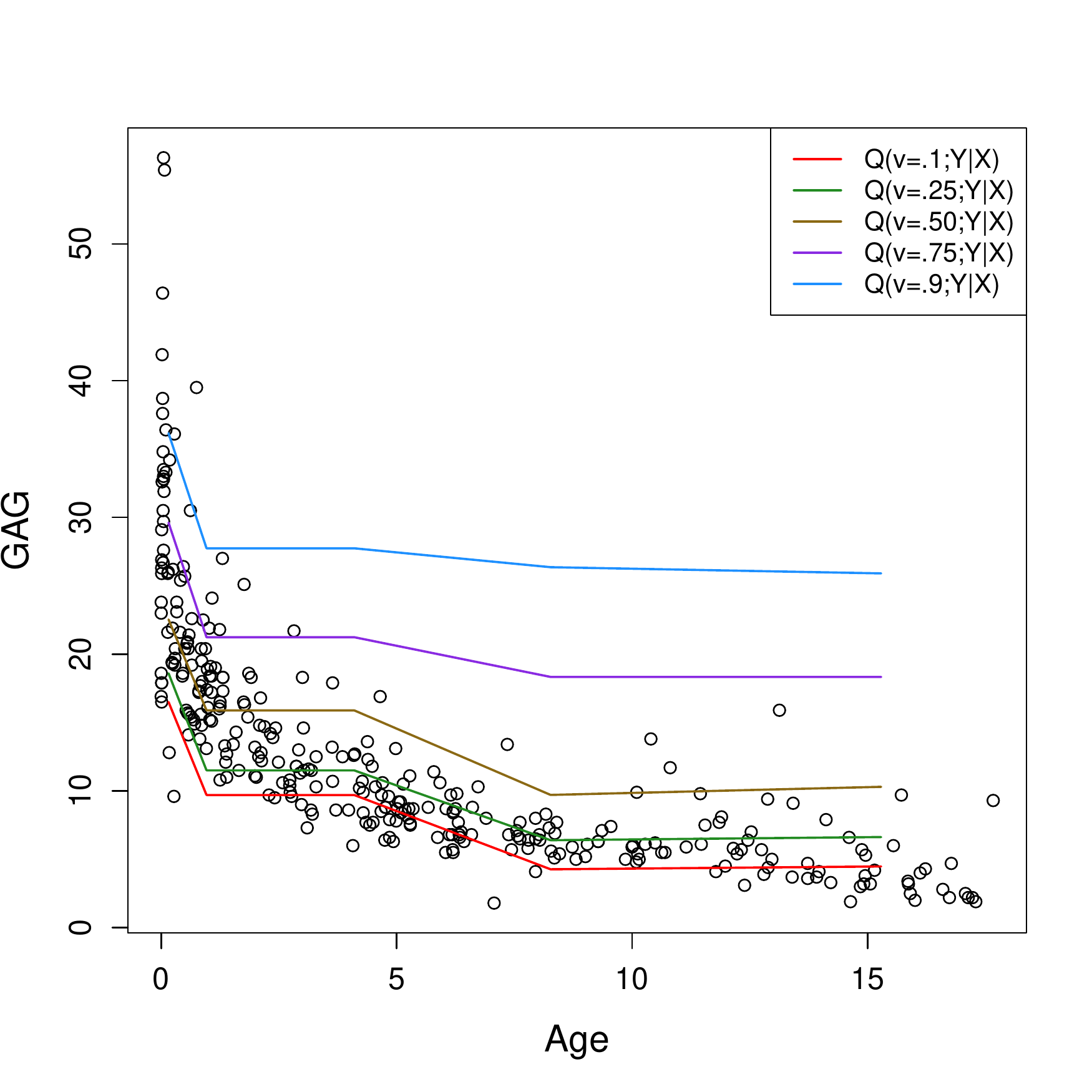}
   }

 \caption[]{(color online) Smooth nonparametric copula based nonlinear regression. Conditional quantile  }\label{fig:gagreg}
\end{figure}

\begin{figure}[tbh!]
 \centering
 \includegraphics[height=\textheight,width=.7\textwidth,keepaspectratio,trim=.5cm .5cm .5cm .5cm]{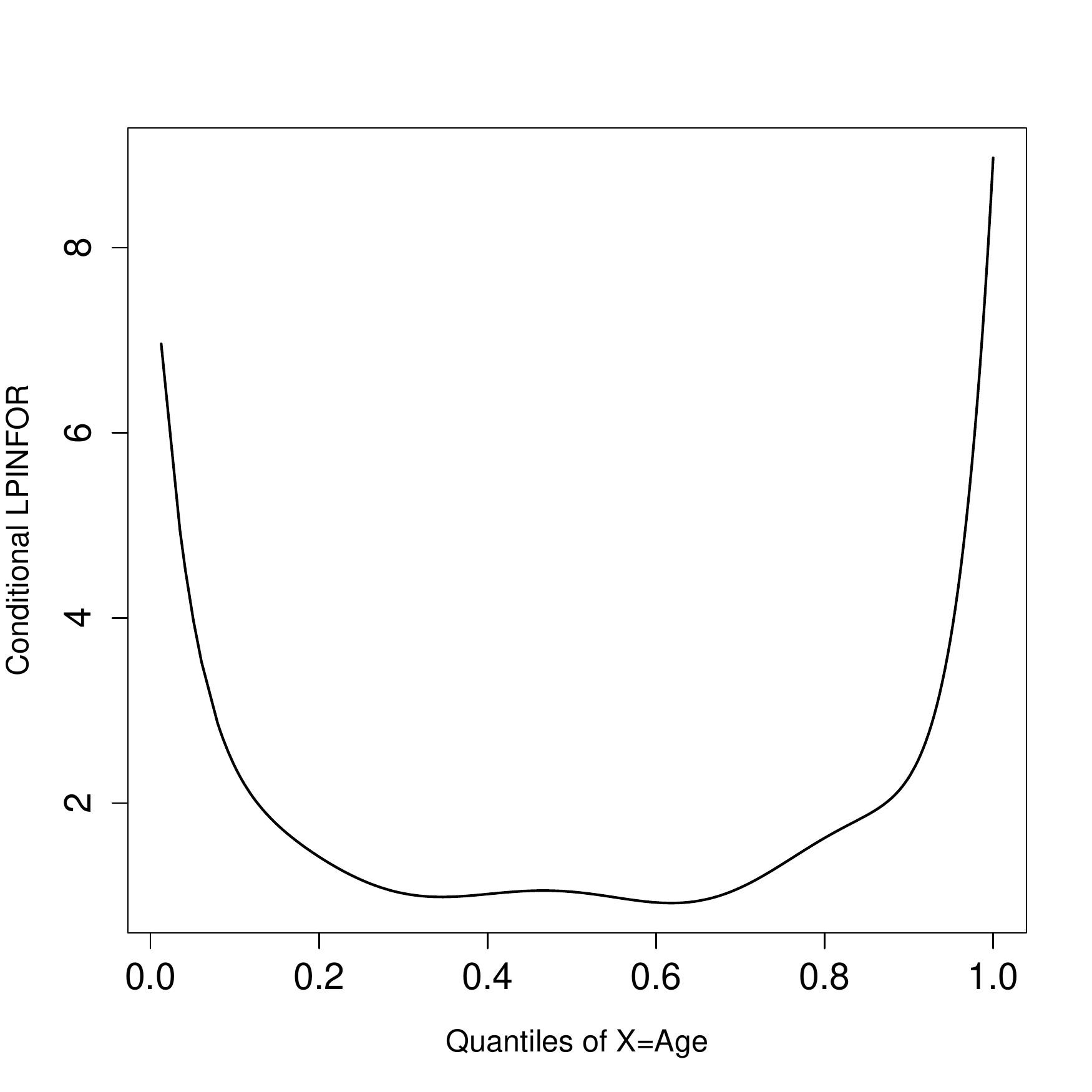}
 \includegraphics[height=\textheight,width=.7\textwidth,keepaspectratio,trim=.5cm .5cm .5cm .5cm]{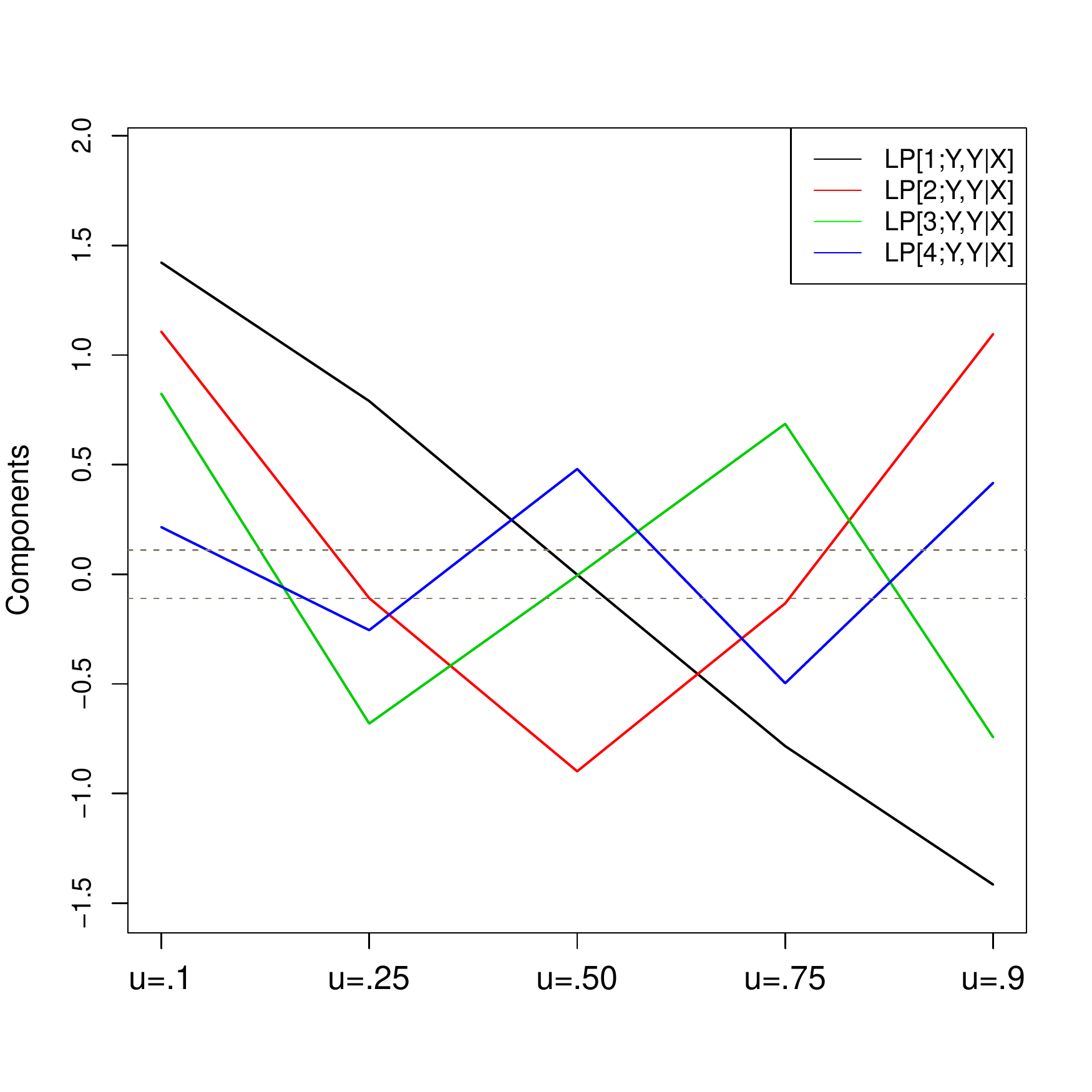}
\vspace{-.1em}
\caption{(color online) Conditional LPINFOR curve and its component decomposition for Ripley data. }
\label{fig:clpi}
\end{figure}


\begin{figure}[tbh!]
 \centering
 \includegraphics[height=\textheight,width=\textwidth,keepaspectratio,trim=1.5cm .5cm 0cm 1.5cm]{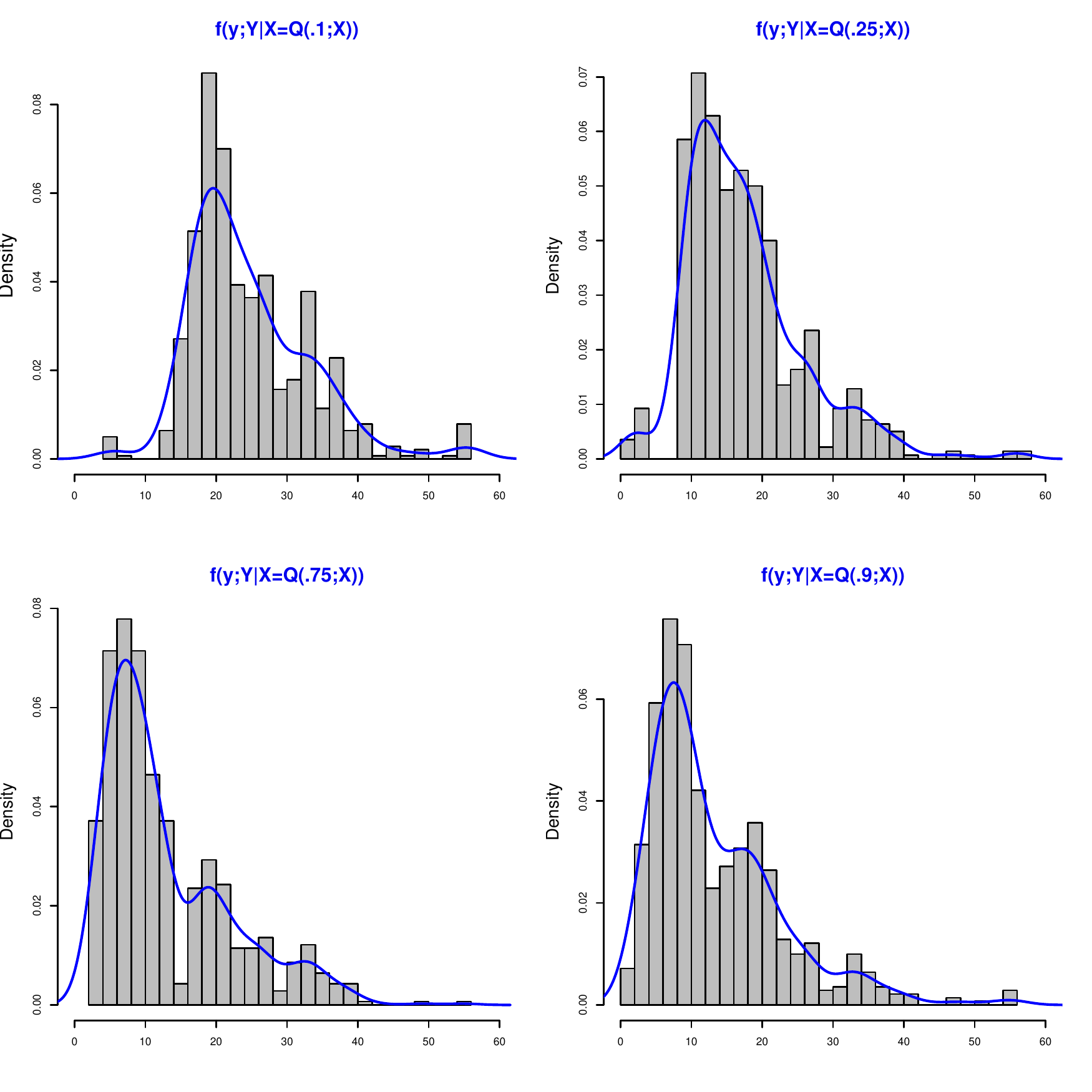}
\vspace{-.1em}
\caption{The estimated nonparametric conditional distributions $f(y;Y|X=Q(u;X))$ for $u=.1,.25,.75,.9$.} \label{fig:gagcd}
\end{figure}

\subsection{Conditional Distribution, Conditional Quantile}
Given $(X,Y)$, it is often of interest to estimate the conditional quantiles and conditional density for comprehensive insight. Our nonparametric estimation of conditional density $f(y;Y|X=x)$ is based on the following LP Skew modeling for $x,y \in \cR$
\beq
f(y;Y|X=x)\,=\,f(y;Y) \cop(F(x;X),F(y;Y))\,=\, f(y;Y)\, d\big(F(y;Y);Y,Y|X=x\big).
\eeq
LP copula estimation theory (presented in Section 4) allows simultaneous estimation of all the copula slices or conditional comparison densities $d(v;Y|X=Q(u;X))$ for various $u$.

Conditional quantiles $Q(v;Y,Y|X=Q(u;X))$ can be simulated  by accept-reject sampling from conditional comparison density. Our approach generalizes the seminal work by \cite{koenker1978} on parametric (linear) quantile regression to nonparametric nonlinear setup, which guaranteed to produce non-crossing quantile curves (a longstanding practical problem, see \cite{koenker2005,he1997} for details).
\subsection{Example}
A comprehensive analysis of the Ripley data is presented using LP algorithmic modeling. Our goal is to answer the scientific question of Ripley data based on $n = 314$ sample,: what are \emph{normal levels} of GAG in children of each age between $1$ and $18$ years of age. In the following, we describe the steps of LP nonparametric algorithm.

\vskip.25em

{\bf Step A.} \cite{ripley04} discusses various traditional regression modeling approaches using polynomial, spline, local polynomials. Figure \ref{fig:gagreg}(a,b) displays estimated LP nonparametric (copula based) conditional mean $\Ex[Y|X=x]$ function \eqref{eq:cmean} given by:
\[ \widehat{\Ex}[Y|X=x]\,=\, 13.1 - 7.32 \,T_1(x) +2.20\,T_2(x).  \]
\vskip.15em

{\bf Step B.} We simulate the conditional quantiles from the estimated conditional comparison densities (slices of copula density), and the resulting smooth nonlinear conditional quantile curves $Q(v;Y,Y|X=Q(u;X))$ are shown in Fig \ref{fig:gagreg}(c,d)  for $u=.1,.25,.5,.75,.9$. Recall from Fig \ref{fig:marginalgag} that the marginal distributions of $X$ (age) and $Y$ (GAG levels) are highly skewed and heavy tailed, which creates a prominent sparse region at the ``tail'' areas. Our algorithm does a satisfactory job in nonparametrically estimating the extreme quantile curves under this challenging scenario. We argue that these estimated quantile regression curves provide the statistical solution to the scientific question: what are ``normal levels'' of GAG in children of each age between $1$ and $18$?
\vskip.25em
{\bf Step C.} Conditional LPINFOR curve is shown in  Fig \ref{fig:clpi}(a). The decomposition of conditional LPINFOR into the components describes the shape of conditional distribution changes, as shown in  Fig \ref{fig:clpi}(b). For example, $\LP[1;Y,Y|X=Q(u;X)]$ describes the conditional mean, $\LP[2;Y,Y|X=Q(u;X)]$ the conditional variance, and so on. Tail behaviour of the conational distributions seem to change rapidly as captured by the $\LP[4;Y,Y|X=Q(u;X)]$ in Fig \ref{fig:clpi}(b).
\vskip.25em
{\bf Step D.} The estimated conditional distributions $f(y;Y|X=Q(u;X))$ are displayed in Fig \ref{fig:gagcd} for $u=.1,.25,.75,.9$. Conditional density shows interesting bi-modality shape at the extreme quantiles. It is also evident from Figure \ref{fig:gagcd} that the classical location-scale shift regression model is inappropriate for this example, which necessitates going beyond the conditional mean description for modeling this data that can tackle the non-Gaussian heavy-tailed data.
%
%
%
%
%
%
%
%
\vskip2em
\section*{Appendix}

\subsection*{A1. Shifted Orthonormal Legendre Polynomials}
Orthonormal Legendre polynomials $\Leg_j(u)$ on interval $0<u<1$
\beas
\Leg_0(u)&=&1 \\
\Leg_1(u)&=& \sqrt{12} (u-.5)\\
\Leg_1(u)&=& \sqrt{5} (6u^2-6u+1)\\
\Leg_3(u)&=&\sqrt{7} (20u^3-30u^2+12u-1)\\
\Leg_4(u)&=&3(70u^4-140u^3+90u^2-20u+1)\\
&\vdots&
\eeas
\subsection*{A2. Proof of Lemma \ref{lemma:qf1}}
To prove $Q(F(X))=X$ with probability $1$. Call $u$ probable if there exists $x$ such that $F(x;X)=u$. Define $x(u)$ to be smallest $x$ satisfying $F(x;X)=u$. Verify $Q(F(x(u)))=x(u)$ for $u$ probable. With probability $1$ the observed value of $X$ belongs to $\{x:x=x(u)\}$ for some $u$ probable. Therefore with probability $1$, $Q(F(X))=X$.

\subsection*{A3. Table of LP Moments}
Table 5 shows LP moments for some standard discrete and continuous distributions.

\begin{table}[htb!]
\setlength{\tabcolsep}{5pt}
\caption{LP Moments of standard distributions.}
\vskip1.5em
\centering
\begin{tabular}{c c cccccc}
\toprule
&\multicolumn{6}{c}{LP Moments}\\
\cmidrule(r){2-7}
Distribution&$\LP[1;X]$ & $\LP[2;X]$ & $\LP[3;X]$ & $\LP[4;X]$ & $\LP[5;X]$ & $\LP[6;X]$\\
\midrule
Poisson$(\la=2)$ & 1.371 &0.205 &0.225 &0.110 &0.103 &0.073  \\[.24ex]
Geometric$(p=0.2)$ & 3.880 & 1.668 & 0.985 & 0.670 & 0.493 & 0.381\\[.24ex]
Uniform$[0,1]$ &   0.289  &0&0&0&0&0         \\[.24ex]
Standard Normal& 0.977 &0 &0.183 &0 &0.081 &0\\[.24ex]
Student's t-df $2$&  1.926 &0 &1.104 &0 &0.866 &0 \\[.24ex]
Chi-squared-df $4$  &2.598 &0.787 &0.562 &0.324 &0.268 &0.187       \\[.24ex]
\hline
\end{tabular}
\label{tab:time}
\end{table}

\subsection*{A4. Table of LP Comoments}
LP-comoment matrix of order $m=4$ is displayed in the following equation for bivariate standard normal distribution with $\rho=0,.5,.9$.

\[\begin{bmatrix}
  0.0 &0.0 &0.0 &0.0\\
  0.0 &0.0 &0.0 &0.0\\
  0.0 &0.0 &0.0 &0.0\\
  0.0 &0.0 &0.0 &0.0
\end{bmatrix} \hskip2em
\begin{bmatrix}
  0.48 &0.0 &0.07 &0.0\\
  0.0 &0.20 &0.0 &0.07\\
  0.07 &0.0 &0.08 &0.0\\
  0.0 &0.07 &0.0 &0.0
\end{bmatrix} \hskip2em
\begin{bmatrix}
  0.89 &0.0 &0.04 &0.0\\
  0.0 &0.76&0.0 &0.09\\
  0.04 &0.0 &0.61 &0.0\\
  0.0 &0.09 &0.0 &0.47
\end{bmatrix}
\]
\vskip.3em
$\hskip5em \rho=0  \hskip10em \rho=.5 \hskip10em \rho=.9\hskip4.5em$

\vskip2em
\bib


\begin{thebibliography}{}

\bibitem[\protect\citeauthoryear{Benzecri}{Benzecri}{1969}]{benzecri1969}
Benzecri, J.~P. (1969).
\newblock {\em Statistical analysis as a tool to make patterns emerge from
  data}.
\newblock New York: Academic Press.

\bibitem[\protect\citeauthoryear{Best and Rayner}{Best and
  Rayner}{1999}]{best1999}
Best, D. and J.~Rayner (1999).
\newblock Goodness of fit for the poisson distribution.
\newblock {\em Statistics \& probability letters\/}~{\em 44\/}(3), 259--265.

\bibitem[\protect\citeauthoryear{Best and Rayner}{Best and
  Rayner}{2003}]{best2003}
Best, D. and J.~Rayner (2003).
\newblock Tests of fit for the geometric distribution.
\newblock {\em Communications in Statistics-Simulation and Computation\/}~{\em
  32\/}(4), 1065--1078.

\bibitem[\protect\citeauthoryear{Blumentritt and Schmid}{Blumentritt and
  Schmid}{2012}]{blumentritt2012}
Blumentritt, T. and F.~Schmid (2012).
\newblock Nonparametric estimation of copula-based measures of multivariate
  association from contingency tables.
\newblock {\em Journal of Statistical Computation and
  Simulation\/}~(ahead-of-print), 1--17.

\bibitem[\protect\citeauthoryear{Breiman and Friedman}{Breiman and
  Friedman}{1985}]{breiman85}
Breiman, L. and J.~H. Friedman (1985).
\newblock Estimating optimal transformations for multiple regression and
  correlation.
\newblock {\em Journal of the American Statistical Association\/}~{\em
  80\/}(391), 580--598.

\bibitem[\protect\citeauthoryear{Csisz{\'a}r}{Csisz{\'a}r}{1975}]{csiszar1975}
Csisz{\'a}r, I. (1975).
\newblock I-divergence geometry of probability distributions and minimization
  problems.
\newblock {\em The Annals of Probability\/}, 146--158.

\bibitem[\protect\citeauthoryear{D'Agostino}{D'Agostino}{1971}]{dagos1971}
D'Agostino, R.~B. (1971).
\newblock An omnibus test of normality for moderate and large size samples.
\newblock {\em Biometrika\/}~{\em 58\/}(2), 341--348.

\bibitem[\protect\citeauthoryear{D'Agostino}{D'Agostino}{1986}]{d1986}
D'Agostino, R.~B. (1986).
\newblock {\em Goodness-of-fit-techniques}, Volume~68.
\newblock CRC press.

\bibitem[\protect\citeauthoryear{David}{David}{1968}]{david1968}
David, H. (1968).
\newblock Gini's mean difference rediscovered.
\newblock {\em Biometrika\/}~{\em 55\/}(3), 573--575.

\bibitem[\protect\citeauthoryear{Denuit and Lambert}{Denuit and
  Lambert}{2005}]{denuit05}
Denuit, M. and P.~Lambert (2005).
\newblock Constraints on concordance measures in bivariate discrete data.
\newblock {\em Journal of Multivariate Analysis\/}~{\em 93\/}(1), 40--57.

\bibitem[\protect\citeauthoryear{Dette, Van~Hecke, and Volgushev}{Dette
  et~al.}{2013}]{dette2013}
Dette, H., R.~Van~Hecke, and S.~Volgushev (2013).
\newblock Misspecification in copula-based regression.
\newblock {\em arXiv preprint arXiv:1310.8037\/}.

\bibitem[\protect\citeauthoryear{Donoho and Jin}{Donoho and Jin}{2008}]{HC08}
Donoho, D. and J.~Jin (2008).
\newblock Higher criticism thresholding: optimal feature selection when useful
  features are rare and weak.
\newblock {\em Proc. Natl. Acad. Sci. USA,\/}~{\em 105}, 14790--15795.

\bibitem[\protect\citeauthoryear{Downton}{Downton}{1966}]{downton1966}
Downton, F. (1966).
\newblock Linear estimates with polynomial coefficients.
\newblock {\em Biometrika\/}, 129--141.

\bibitem[\protect\citeauthoryear{Eagleson}{Eagleson}{1964}]{eagleson64}
Eagleson, G. (1964).
\newblock Polynomial expansions of bivariate distributions.
\newblock {\em The Annals of Mathematical Statistics\/}~{\em 35\/}(3),
  1208--1215.

\bibitem[\protect\citeauthoryear{Emerson}{Emerson}{1968}]{emerson1968}
Emerson, P.~L. (1968).
\newblock Numerical construction of orthogonal polynomials from a general
  recurrence formula.
\newblock {\em Biometrics\/}, 695--701.

\bibitem[\protect\citeauthoryear{Fisher}{Fisher}{1940}]{fisher1940}
Fisher, R.~A. (1940).
\newblock The precision of discriminant functions.
\newblock {\em Annals of Eugenics\/}~{\em 10\/}(1), 422--429.

\bibitem[\protect\citeauthoryear{Genest and Neslehova}{Genest and
  Neslehova}{2007}]{genest2007}
Genest, C. and J.~Neslehova (2007).
\newblock A primer on copulas for count data.
\newblock {\em Astin Bulletin\/}~{\em 37\/}(2), 475--515.

\bibitem[\protect\citeauthoryear{Genest, Ne{\v{s}}lehov{\'a}, and
  R{\'e}millard}{Genest et~al.}{2013}]{genest2013}
Genest, C., J.~Ne{\v{s}}lehov{\'a}, and B.~R{\'e}millard (2013).
\newblock On the estimation of spearman's rho and related tests of independence
  for possibly discontinuous multivariate data.
\newblock {\em Journal of Multivariate Analysis\/}, 214--228.

\bibitem[\protect\citeauthoryear{Gilula, Krieger, and Ritov}{Gilula
  et~al.}{1988}]{gilula1988}
Gilula, Z., A.~M. Krieger, and Y.~Ritov (1988).
\newblock Ordinal association in contingency tables: some interpretive aspects.
\newblock {\em Journal of the American Statistical Association\/}~{\em
  83\/}(402), 540--545.

\bibitem[\protect\citeauthoryear{Goodman}{Goodman}{1981}]{goodman1981}
Goodman, L.~A. (1981).
\newblock Association models and canonical correlation in the analysis of
  cross-classifications having ordered categories.
\newblock {\em Journal of the American Statistical Association\/}~{\em
  76\/}(374), 320--334.

\bibitem[\protect\citeauthoryear{Goodman}{Goodman}{1985}]{goodman1985}
Goodman, L.~A. (1985).
\newblock The analysis of cross-classified data having ordered and/or unordered
  categories: Association models, correlation models, and asymmetry models for
  contingency tables with or without missing entries.
\newblock {\em The Annals of Statistics\/}~{\em 13}, 10--69.

\bibitem[\protect\citeauthoryear{Goodman}{Goodman}{1991}]{goodman1991}
Goodman, L.~A. (1991).
\newblock Measures, models, and graphical displays in the analysis of
  cross-classified data (with discussion).
\newblock {\em Journal of the American Statistical association\/}~{\em
  86\/}(416), 1085--1111.

\bibitem[\protect\citeauthoryear{Goodman}{Goodman}{1996}]{good96}
Goodman, L.~A. (1996).
\newblock A single general method for the analysis of cross-classified data:
  reconciliation and synthesis of some methods of {P}earson, {Y}ule, and
  {F}isher, and also some methods of correspondence analysis and association
  analysis.
\newblock {\em Journal of the American Statistical Association\/}~{\em
  91\/}(433), 408--428.

\bibitem[\protect\citeauthoryear{Greenacre}{Greenacre}{2010}]{greenacre2010}
Greenacre, M. (2010).
\newblock {\em Correspondence analysis in practice}.
\newblock CRC Press.

\bibitem[\protect\citeauthoryear{Haberman}{Haberman}{1988}]{haberman88}
Haberman, S.~J. (1988).
\newblock A warning on the use of chi-squared statistics with frequency tables
  with small expected cell counts.
\newblock {\em Journal of the American Statistical Association\/}~{\em
  83\/}(402), 555--560.

\bibitem[\protect\citeauthoryear{He}{He}{1997}]{he1997}
He, X. (1997).
\newblock Quantile curves without crossing.
\newblock {\em The American Statistician\/}~{\em 51\/}(2), 186--192.

\bibitem[\protect\citeauthoryear{Hirschfeld}{Hirschfeld}{1935}]{hirschfeld1935}
Hirschfeld, H.~O. (1935).
\newblock A connection between correlation and contingency.
\newblock {\em Proceedings of the Cambridge Philosophical Society\/}~{\em 31},
  520--524.

\bibitem[\protect\citeauthoryear{Hosking}{Hosking}{1990}]{hos90}
Hosking, J.~R. (1990).
\newblock L-moments: analysis and estimation of distributions using linear
  combinations of order statistics.
\newblock {\em Journal of the Royal Statistical Society. Series B
  (Methodological)\/}, 105--124.

\bibitem[\protect\citeauthoryear{Hosking and Wallis}{Hosking and
  Wallis}{1997}]{HosL}
Hosking, J. R.~M. and J.~R. Wallis (1997).
\newblock {\em Regional Frequency Analysis: An Approach Based on L-moments.}
\newblock Cambridge University Press, Cambridge.

\bibitem[\protect\citeauthoryear{Kallenberg and Ledwina}{Kallenberg and
  Ledwina}{1999}]{kall99}
Kallenberg, W.~C. and T.~Ledwina (1999).
\newblock Data-driven rank tests for independence.
\newblock {\em Journal of the American Statistical Association\/}~{\em
  94\/}(445), 285--301.

\bibitem[\protect\citeauthoryear{Khmaladze}{Khmaladze}{2013}]{khmaladze2013}
Khmaladze, E. (2013).
\newblock Note on distribution free testing for discrete distributions.
\newblock {\em The Annals of Statistics\/}~{\em 41\/}(6), 2979--2993.

\bibitem[\protect\citeauthoryear{Koenker}{Koenker}{2005}]{koenker2005}
Koenker, R. (2005).
\newblock {\em Quantile regression}.
\newblock Number~38. Cambridge university press.

\bibitem[\protect\citeauthoryear{Koenker and Bassett}{Koenker and
  Bassett}{1978}]{koenker1978}
Koenker, R. and G.~Bassett (1978).
\newblock Regression quantiles.
\newblock {\em Econometrica: journal of the Econometric Society\/}, 33--50.

\bibitem[\protect\citeauthoryear{Koziol}{Koziol}{1979}]{koziol1979}
Koziol, J.~A. (1979).
\newblock A smooth test for bivariate independence.
\newblock {\em Sankhy{\=a}: The Indian Journal of Statistics, Series B\/},
  260--269.

\bibitem[\protect\citeauthoryear{Lewi}{Lewi}{1998}]{lewi1998}
Lewi, P. (1998).
\newblock Analysis of contingency tables.
\newblock {\em Handbook of Chemometrics and Qualimetrics\/}~(Part B), 161--206.

\bibitem[\protect\citeauthoryear{Mack and Wolfe}{Mack and
  Wolfe}{1981}]{mack1981}
Mack, G.~A. and D.~A. Wolfe (1981).
\newblock K-sample rank tests for umbrella alternatives.
\newblock {\em Journal of the American Statistical Association\/}~{\em
  76\/}(373), 175--181.

\bibitem[\protect\citeauthoryear{Mesfioui and Quessy}{Mesfioui and
  Quessy}{2010}]{mesfioui2010}
Mesfioui, M. and J.-F. Quessy (2010).
\newblock Concordance measures for multivariate non-continuous random vectors.
\newblock {\em Journal of Multivariate Analysis\/}~{\em 101\/}(10), 2398--2410.

\bibitem[\protect\citeauthoryear{Mukhopadhyay}{Mukhopadhyay}{2013}]{D11d}
Mukhopadhyay, S. (2013).
\newblock Cdfdr: A comparison density approach to local false discovery rate
  estimation.
\newblock {\em arXiv:1308.2403\/}.

\bibitem[\protect\citeauthoryear{Mukhopadhyay and Parzen}{Mukhopadhyay and
  Parzen}{2013}]{D12e}
Mukhopadhyay, S. and E.~Parzen (2013).
\newblock Nonlinear time series modeling by {LPT}ime, nonparametric empirical
  learning.
\newblock {\em arXiv:1308.0642\/}.

\bibitem[\protect\citeauthoryear{Ne{\v{s}}lehov{\'a}}{Ne{\v{s}}lehov{\'a}}{2007}]{nevslehova07}
Ne{\v{s}}lehov{\'a}, J. (2007).
\newblock On rank correlation measures for non-continuous random variables.
\newblock {\em Journal of Multivariate Analysis\/}~{\em 98\/}(3), 544--567.

\bibitem[\protect\citeauthoryear{Neyman}{Neyman}{1937}]{Neyman37}
Neyman, J. (1937).
\newblock Smooth tests for goodness of fit.
\newblock {\em Skand. Aktuar.\/}~{\em 20}, 150--199.

\bibitem[\protect\citeauthoryear{Nikoloulopoulos and Karlis}{Nikoloulopoulos
  and Karlis}{2010}]{nikolo2010}
Nikoloulopoulos, A.~K. and D.~Karlis (2010).
\newblock Regression in a copula model for bivariate count data.
\newblock {\em Journal of Applied Statistics\/}~{\em 37\/}(9), 1555--1568.

\bibitem[\protect\citeauthoryear{Noah and Tibshirani}{Noah and
  Tibshirani}{2011}]{Tibs11}
Noah, S. and R.~Tibshirani (2011).
\newblock Comment on detecting novel associations in large data sets.
\newblock {\em Science Comment\/}.

\bibitem[\protect\citeauthoryear{Noh, Ghouch, and Bouezmarni}{Noh
  et~al.}{2013}]{noh2013}
Noh, H., A.~E. Ghouch, and T.~Bouezmarni (2013).
\newblock Copula-based regression estimation and inference.
\newblock {\em Journal of the American Statistical Association\/}~{\em
  108\/}(502), 676--688.

\bibitem[\protect\citeauthoryear{Parzen}{Parzen}{1962}]{parzen1962}
Parzen, E. (1962).
\newblock On estimation of a probability density function and mode.
\newblock {\em Annals of mathematical statistics\/}~{\em 33\/}(3), 1065--1076.

\bibitem[\protect\citeauthoryear{Parzen}{Parzen}{1979}]{parzen79}
Parzen, E. (1979).
\newblock Nonparametric statistical data modeling (with discussion).
\newblock {\em Journal of the American Statistical Association\/}~{\em 74},
  105--131.

\bibitem[\protect\citeauthoryear{Parzen}{Parzen}{2004}]{parzen04b}
Parzen, E. (2004).
\newblock Quantile probability and statistical data modeling.
\newblock {\em Statistical Science,\/}~{\em 19}, 652--662.

\bibitem[\protect\citeauthoryear{Pearson}{Pearson}{1900}]{pearson1900}
Pearson, K. (1900).
\newblock On the criterion that a given system of deviations from the probable
  in the case of a correlated system of variables is such that it can be
  reasonably supposed to have arisen from random sampling.
\newblock {\em The London, Edinburgh, and Dublin Philosophical Magazine and
  Journal of Science\/}~{\em 50\/}(302), 157--175.

\bibitem[\protect\citeauthoryear{Rayner and Best}{Rayner and
  Best}{1996}]{rayner96}
Rayner, J. and D.~Best (1996).
\newblock Smooth extensions of {P}earsons's product moment correlation and
  spearman's rho.
\newblock {\em Statistics \& Probability Letters\/}~{\em 30\/}(2), 171--177.

\bibitem[\protect\citeauthoryear{Reshef, Reshef, Finucane, Grossman, McVean,
  Turnbaugh, Lander, and Mitzenmacher}{Reshef et~al.}{2011}]{MI11}
Reshef, D.~N., Y.~A. Reshef, H.~K. Finucane, S.~R. Grossman, G.~McVean, P.~J.
  Turnbaugh, E.~S. Lander, and M.~Mitzenmacher (2011).
\newblock Detecting novel associations in large data sets.
\newblock {\em Science\/}~{\em 334\/}(6062), 1518--1524.

\bibitem[\protect\citeauthoryear{Ripley}{Ripley}{2004}]{ripley04}
Ripley, B.~D. (2004).
\newblock Selecting amongst large classes of models.
\newblock {\em URL http://www.stats.ox.ac.uk/~ripley/Nelder80.pdf. Symposium in
  Honour of David Cox's 80th birthday\/}.

\bibitem[\protect\citeauthoryear{Schechtman and Yitzhaki}{Schechtman and
  Yitzhaki}{1987}]{gini87}
Schechtman, E. and S.~Yitzhaki (1987).
\newblock A measure of association based on {G}ini's mean difference.
\newblock {\em Communications in Statistics - Theory and Methods\/}~{\em
  A16\/}(1), 207--231.

\bibitem[\protect\citeauthoryear{Serfling and Xiao}{Serfling and
  Xiao}{2007}]{Serf07}
Serfling, R. and P.~Xiao (2007).
\newblock A contribution to multivariate {L}-moments: L-comoment matrices.
\newblock {\em Journal of Multivariate Analysis\/}~{\em 98\/}(9), 1765--1781.

\bibitem[\protect\citeauthoryear{Sillitto}{Sillitto}{1969}]{sillitto1969}
Sillitto, G. (1969).
\newblock Derivation of approximants to the inverse distribution function of a
  continuous univariate population from the order statistics of a sample.
\newblock {\em Biometrika\/}~{\em 56\/}(3), 641--650.

\bibitem[\protect\citeauthoryear{Sklar}{Sklar}{1959}]{sklar1959}
Sklar, M. (1959).
\newblock Fonctions de r{\'e}partition {\`a} n dimensions et leurs marges.
\newblock {\em Publ. Inst. Statistique Univ. Paris\/}~{\em 8}, 229--231.

\bibitem[\protect\citeauthoryear{Sur, Shmueli, Bose, and Dubey}{Sur
  et~al.}{2013}]{Shmueli2013}
Sur, P., G.~Shmueli, S.~Bose, and P.~Dubey (2013).
\newblock Modeling bimodal discrete data using {C}onway-{M}axwell-{P}oisson
  mixture models.
\newblock {\em arXiv:1309.0579\/}.

\bibitem[\protect\citeauthoryear{Sz{\'e}kely and Rizzo}{Sz{\'e}kely and
  Rizzo}{2009}]{szekely09}
Sz{\'e}kely, G.~J. and M.~L. Rizzo (2009).
\newblock Brownian distance covariance.
\newblock {\em The annals of applied statistics\/}, 1236--1265.

\bibitem[\protect\citeauthoryear{Yitzhaki and Schechtman}{Yitzhaki and
  Schechtman}{2013}]{yitzhaki13book}
Yitzhaki, S. and E.~Schechtman (2013).
\newblock {\em The {G}ini Methodology: A primer on a statistical methodology},
  Volume 272.
\newblock Springer.

\end{thebibliography}
\end{document}